# Monotone Control Systems

David Angeli
Dip. Sistemi e Informatica
University of Florence, 50139 Firenze, Italy
angeli@dsi.unifi.it

Eduardo D. Sontag[*]
Dept. of Mathematics
Rutgers University, NJ, USA
sontag@hilbert.rutgers.edu

**Abstract**

Monotone systems constitute one of the most important classes of dynamical systems used in mathematical biology modeling. The objective of this paper is to extend the notion of monotonicity to systems with inputs and outputs, a necessary first step in trying to understand interconnections, especially including feedback loops, built up out of monotone components. Basic definitions and theorems are provided, as well as an application to the study of a model of one of the cell's most important subsystems.

## 1 Introduction

One of the most important classes of dynamical systems in theoretical biology is that of *monotone systems*. Among the classical references in this area are the textbook by Smith [24] and the papers [14, 15] by Hirsh and [23] by Smale. Monotone systems are those for which trajectories preserve a partial ordering on states. They include the subclass of *cooperative* systems (see e.g. [1, 6, 7] for recent contributions in the control literature), for which different state variables reinforce each other (positive feedback) as well as more general systems in which each pair of variables may affect each other in either positive or negative, or even mixed, forms (precise definitions are given below). Although one may consider systems in which constant parameters (which can be thought of as constant inputs) appear, as done in the recent paper [22] for cooperative systems, the concept of monotone system has been traditionally defined only for systems with no external input (or "control") *functions*.

The objective of this paper is to extend the notion of monotone systems to *systems with inputs and outputs*. This is by no means a purely academic exercise, but it is a necessary first step in trying to understand interconnections, especially including feedback loops, built up out of monotone components.

The successes of systems theory have been due in large part to its ability to analyze complicated structures on the basis of the behavior of elementary subsystems, each of which is "nice" in a suitable input/output sense (stable, passive, etc), in conjunction with the use of tools such as the small gain theorem to characterize interconnections.

On the other hand, one of the main themes and challenges in current molecular biology lies in the understanding of cell behavior in terms of cascade and feedback interconnections of elementary "modules" which appear repeatedly, see e.g. [13]. Our work reported here was motivated by the problem of studying one such module type (closely related to, but more general than, the example which motivated [27]), and the realization that the theory of monotone systems, when extended to allow for inputs, provides an appropriate tool to formulate and prove basic properties of such modules.

---

[*]Supported in part by US Air Force Grant F49620-01-1-0063, by a grant from the National Institutes of Health (P20 GM64375) and by funding from Rutgers University to the Institute for Biology at the Interface of the Mathematical and Physical Sciences (BioMaPS Institute).

The organization of this paper is as follows. In Section 2, we introduce the basic concepts, including the special case of cooperative systems. Section 3 provides infinitesimal characterizations of monotonicity, relying upon certain technical points discussed in the Appendix. Cascades are the focus of Section 4, and Section 5 introduces the notions of static Input/State and Input/Output characteristics, which then play a central role in the study of feedback interconnections and a small-gain theorem in Section 6. We return to the biological example of MAPK cascades in Section 7. Finally, Section 8 shows the equivalence between cooperative systems and positivity of linearizations.

We view this paper as only the begining of a what should be a fruitful direction of research into a new type of nonlinear systems. In particular, in a follow-up paper [3], we present results dealing with more general feedback interconnections than those considered in this paper, as well as graphical criteria for monotonicity.

## 2 Monotone Systems

Monotone dynamical systems are usually defined on subsets of ordered Banach (or even more general metric) spaces. An *ordered Banach space* is a real Banach space $\mathbb{B}$ together with a distinguished nonempty closed subset $K$ of $\mathbb{B}$, its *positive cone*. (The spaces $\mathbb{B}$ which we study in this paper will all be Euclidean spaces; however, the basic definitions can be given in more generality, and doing so might eventually be useful for applications such as the study of systems with delays, as done in [24] for systems without inputs.) The set $K$ is assumed to have the following properties: it is a cone, i.e. $\alpha K \subset K$ for $\alpha \in \mathbb{R}_+$, it is convex (equivalently, since $K$ is a cone, $K + K \subset K$), and $K \cap (-K) = \{0\}$. An ordering is then defined by $x_1 \succeq x_2$ iff $x_1 - x_2 \in K$. Strict ordering is denoted by $x_1 \succ x_2$, meaning that $x_1 \succeq x_2$ and $x_1 \neq x_2$. One often uses as well the notations $\prec$ and $\preceq$, in the obvious sense ($x_2 \preceq x_1$ means $x_1 \succeq x_2$). (Most of the results discussed in this paper use only that $K$ is a cone. The property $K \cap (-K) = \{0\}$, which translates into reflexivity of the order, is never used, and the convexity property, which translates into transitivity of the order, will be only used in a few places.)

The most typical example would be $\mathbb{B} = \mathbb{R}^n$ and $K = \mathbb{R}^n_{\geq 0}$, in which case "$x_1 \succeq x_2$" means that each coordinate of $x_1$ is bigger or equal than the corresponding coordinate of $x_2$. This order on state spaces gives rise to the class of "cooperative systems" discussed below. However, other orthants in $\mathbb{R}^n$ other than the positive orthant $K = \mathbb{R}^n_{\geq 0}$ are often more natural in applications, as we will see.

In view of our interest in biological and chemical applications, we must allow state spaces to be non-linear subsets of linear spaces. For example, state variables typically represent concentrations, and hence must be positive, and are often subject to additional inequality constraints such as stoichiometry or mass preservation. Thus, from now on, we will assume given an ordered Banach space $\mathbb{B}$ and a subset $X$ of $\mathbb{B}$ which is the closure of an open subset of $\mathbb{B}$. For instance, $X = \mathbb{B}$, or, in an example to be considered later, $\mathbb{B} = \mathbb{R}^2$ with the order induced by $K = \mathbb{R}_{\geq 0} \times \mathbb{R}_{\leq 0}$, and $X = \{(x, y) \in \mathbb{R}^2 \mid x \geq 0, y \geq 0, x + y \leq 1\}$.

We start by recalling the standard concept of monotonicity for uncontrolled systems.

**Definition 2.1** A dynamical system $\phi : \mathbb{R}_{\geq 0} \times X \to X$ is *monotone* if the implication below holds:

$$x_1 \succeq x_2 \Rightarrow \phi(t, x_1) \succeq \phi(t, x_2) \qquad \forall t \geq 0\,. \qquad \square$$

If the interior of the positive cone $K$ is nonempty (as is often the case in applications of monotonicity), one can also define a stricter ordering: $x_1 \gg x_2 \Leftrightarrow x - y \in \text{int}(K)$. (For example, when $K = \mathbb{R}^n_{\geq 0}$, this means that every coordinate of $x_1$ is strictly larger than the corresponding coordinate of $x_2$, in contrast to "$x_1 \succ x_2$" which means merely that some coordinate is strictly bigger while the rest are bigger or equal.) Accordingly, one says that a dynamical system $\phi : \mathbb{R}_{\geq 0} \times X \to X$ is *strictly monotone* if $x_1 \succ x_2$ implies that $\phi(t, x_1) \gg \phi(t, x_2)$ for all $t \geq 0$.

**Inputs**

Next we generalize, in a very natural way, the above definition to *controlled* dynamical systems, i.e., systems forced by some exogenous input signal. In order to do so, we assume given a partially ordered input value space $\mathcal{U}$. Technically, we will assume that $\mathcal{U}$ is a subset of an ordered Banach space $\mathbb{B}_\mathcal{U}$. Thus, for any pair of input values $u_1$ and $u_2 \in \mathcal{U}$, we write $u_1 \succeq u_2$ whenever $u_1 - u_2 \in K_u$ where $K_u$ is the corresponding positivity cone in $\mathbb{B}_\mathcal{U}$.

By an "input" or "control" we shall mean a Lebesgue measurable function $u(\cdot) : \mathbb{R}_{\geq 0} \to \mathcal{U}$ which is essentially bounded, i.e. there is for each finite interval $[0, T]$ some compact subset $C \subseteq \mathcal{U}$ such that $u(t) \in C$ for almost all $t \in [0, T]$. We denote by $\mathcal{U}_\infty$ the set of all inputs. Accordingly, given two $u_1, u_2 \in \mathcal{U}_\infty$, we write[1]
$$u_1 \succeq u_2 \iff u_1(t) \succeq u_2(t) \ \forall \, t \geq 0 \,.$$
A *controlled dynamical system* is specified by a state space $X$ as above, an input set $\mathcal{U}$, and a mapping $\phi : \mathbb{R}_{\geq 0} \times X \times \mathcal{U}_\infty \to X$ such that the usual semigroup properties hold. (Namely, $\phi(0, x, u) = x$ and $\phi(t, \phi(s, x, u_1), u_2) = \phi(s + t, x, v)$, where $v$ is the restriction of $u_1$ to the interval $[0, s]$ concatenated with the restriction of $u_2$ to $[s, \infty)$; we will soon specialize to solutions of controlled differential equations.)

We interpret $\phi(t, \xi, u)$ as the state at time $t$ obtained if the initial state is $\xi$ and the external input is $u(\cdot)$. Sometimes, when clear from the context, we write "$x(t, \xi, u)$" or just "$x(t)$" instead of $\phi(t, \xi, u)$. When there is no risk of confusion, we use "$x$" to denote states (i.e., elements of $X$) as well as trajectories, but for emphasis we sometimes use $\xi$, possibly subscripted, and other Greek letters, to denote states. Similarly, "$u$" may refer to an input value (element of $\mathcal{U}$) or an input function (element of $\mathcal{U}_\infty$).

**Definition 2.2** A controlled dynamical system $\phi : \mathbb{R}_{\geq 0} \times X \times \mathcal{U}_\infty \to X$ is *monotone* if the implication below holds:
$$u_1 \succeq u_2, \ x_1 \succeq x_2 \implies \phi(t, x_1, u_1) \succeq \phi(t, x_2, u_2) \ \forall \, t \geq 0 \,. \qquad \square$$

Viewing systems with no inputs as controlled systems for which the input value space $\mathcal{U}$ has just one element, one recovers the classical definition.

One may also, of course, define strictly monotone controlled dynamical systems, by analogy with the concept defined for systems with no inputs.

### Outputs

We will also consider monotone systems *with outputs* $y = h(x)$. These are specified by a controlled monotone system $\phi$ together with a monotone ($x_1 \succeq x_2 \Rightarrow h(x_1) \succeq h(x_2)$) map $h : X \to \mathcal{Y}$, where $\mathcal{Y}$, the set of measurement or output values, is a subset of some ordered Banach space $\mathbb{B}_\mathcal{Y}$. We often use the shorthand $y(t, x, u)$ instead of $h(\phi(t, x, u))$, to denote the output at time $t$ corresponding to the state obtained from initial state $x$ and input $u$.

### Systems without inputs

For each constant input signal $u(t) \equiv \bar{u}$, a controlled monotone dynamical system behaves as a classical monotone dynamical system. This allows application of the rich theory developed for this class of systems. For example, for strictly monotone dynamical systems, the following properties are satisfied (see [24] for precise statements and proofs):

- For $a, b \in X$, with $b \succ a$, let $[a, b] = \{x \in X : b \succeq x \succeq a\}$. If $a$ and $b$ are equilibria for a strongly monotone system $\phi$ and $\phi([a, b], t)$ is relatively compact for all $t > 0$ then one of following three alternatives must hold:
  (i) all trajectories with $b \succ x \succeq a$ converge to $a$;
  (ii) all trajectories with $b \succeq x \succ a$ converge to $b$; or
  (iii) there exists an equilibrium $c$ with $a \succ c \succ b$

  (order interval trychotomy, [20, 24]).

- If all semiorbits have compact closures, then convergence (as $t \to +\infty$) to the set of equilibria is in some sense "generic" (see [14]).

---
[1]To be more precise, this and other definitions should be interpreted in an "almost everywhere" sense, since inputs are Lebesgue-measurable functions.

- Periodic orbits cannot be stable.
- A theorem of Smale shows that arbitrarily complex dynamics can arise on unstable $n-1$ dimensional invariant subsets of monotone systems.

We will study generalizations of some of these properties to controlled monotone systems in forthcoming work.

### Systems defined by differential equations

From now on, we will specialize to the case of systems defined by differential equations with inputs:

$$\dot{x} = f(x, u) \tag{1}$$

(see [26] for basic definitions and properties regarding such systems). We make the following technical assumptions. The map $f$ is defined on $\widetilde{X} \times \mathcal{U}$, where $\widetilde{X}$ is some open subset of $\mathbb{B}$ which contains $X$, and $\mathbb{B} = \mathbb{R}^n$ for some integer $n$. We assume that $f(x, u)$ is continuous in $(x, u)$ and locally Lipschitz continuous in $x$ locally uniformly on $u$. This last property means that for each compact subsets $C_1 \subseteq X$ and $C_2 \subseteq \mathcal{U}$ there exists some constant $k$ such that $|f(\xi, u) - f(\zeta, u)| \leq |\xi - \zeta|$ for all $\xi, \zeta \in C_1$ and all $u \in C_2$. In order to obtain a well-defined controlled dynamical system on $X$, we will assume that the solution $x(t) = \phi(t, x_0, u)$ of $\dot{x} = f(x, u)$ with initial condition $x(0) = x_0$ is defined for all inputs $u(\cdot)$ and all times $t \geq 0$. This means that solutions with initial states in $X$ must be defined for all $t \geq 0$ (forward completeness) and that the set $X$ is forward invariant. (Forward invariance of $X$ may be checked using tangent cones at the boundary of $X$, see the Appendix.)

*From now on, all systems will be assumed to be of this form.*

## 3  Infinitesimal Characterizations

For systems (1) defined by controlled differential equations, we will provide an infinitesimal characterization of monotonicity, expressed directly in terms of the vector field, which does not require the explicit computation of solutions. Our result will generalize the well-known Kamke conditions, discussed in [24], Chapter 3. We denote $\mathcal{V} := \text{int } X$, the interior of $X$ (recall that $X$ is the closure $\mathcal{V}$) and impose the following *approximability property* (see [24], Remark 3.1.4):

> for all $\xi_1, \xi_2 \in X$ such that $\xi_1 \succeq \xi_2$, there exist sequences $\{\xi_1^i\}, \{\xi_2^i\} \subseteq \mathcal{V}$
> such that $\xi_1^i \succeq \xi_2^i$ for all $i$ and $\xi_1^i \to \xi_1$ and $\xi_2^i \to \xi_2$ as $i \to \infty$.

**Remark 3.1** The approximability assumption is very mild. It is satisfied, in particular, if the set $X$ is convex, and, even more generally, if it is strictly star-shaped with respect to some interior point $\xi_*$, i.e., for all $\xi \in X$ and all $0 \leq \lambda < 1$, it holds that $\lambda \xi + (1-\lambda)\xi_* \in \mathcal{V}$. (Convex sets with nonempty interior have this property with respect to any point $\xi_* \in \mathcal{V}$, since $\lambda \xi + (1-\lambda)\xi_* \in Q := \lambda \xi + (1-\lambda)\mathcal{V} \subseteq X$ (the inclusion by convexity) and the set $Q$ is open because $\eta \mapsto \lambda \xi + (1-\lambda)\eta$ is an invertible affine mapping.) Indeed, suppose that $\xi_1 - \xi_2 \in K$, pick any sequence $\lambda^i \nearrow 1$, and define $\xi_j^i := \lambda^i \xi_j + (1-\lambda^i)\xi_*$ for $j = 1, 2$. These elements are in $\mathcal{V}$, they converge to $\xi_1$ and $\xi_2$ respectively, and each $\xi_1^i - \xi_2^i = \lambda^i(\xi_1 - \xi_2)$ belongs to $K$ because $K$ is a cone. Moreover, a slightly stronger property holds as well, for starshaped $X$, namely: if $\xi_1, \xi_2 \in X$ are such that $\xi_1 \succeq \xi_2$ and if for some linear map $L : \mathbb{R}^n \to \mathbb{R}^q$ it holds that $L\xi_1 = L\xi_2$, then the sequences $\{\xi_1^i\}, \{\xi_2^i\}$ can be picked such that $L\xi_1^i = L\xi_2^i$ for all $i$; this follows from the construction, since $L(\xi_1^i - \xi_2^i) = \lambda^i L(\xi_1 - \xi_2) = 0$. For instance, $L$ might select those coordinates which belong in some subset $I \subseteq \{1, \ldots, n\}$. This stronger property will be useful later, when we look at boundary points. □

The characterization will be in terms of a standard notion of tangent cone, studied in nonsmooth analysis, which we introduce next.

**Definition 3.2** Let $S$ be a subset of a Euclidean space, and pick any $\xi \in \mathcal{S}$. The *tangent cone to* $\mathcal{S}$ *at* $\xi$ is the set

$$\mathcal{T}_\xi \mathcal{S} := \left\{ \lim_{i \to \infty} \frac{1}{t_i}(\xi_i - \xi) \,\middle|\, \xi_i \underset{\mathcal{S}}{\to} \xi, \, t_i \searrow 0 \right\}$$

where "$\xi_i \underset{S}{\to} \xi$" means that $\xi_i \to \xi$ as $i \to \infty$ and that $\xi_i \in S$ for all $i$. □

Several properties of tangent cones are reviewed in the Appendix. The main result in this section is as follows.

**Theorem 1** *The system (1) is monotone if and only if, for all $\xi_1, \xi_2 \in \mathcal{V}$:*

$$\xi_1 \succeq \xi_2 \text{ and } u_1 \succeq u_2 \;\Rightarrow\; f(\xi_1, u_1) - f(\xi_2, u_2) \in \mathcal{T}_{\xi_1 - \xi_2} K \tag{2}$$

*or, equivalently,*

$$\xi_1 - \xi_2 \in \partial K \text{ and } u_1 \succeq u_2 \;\Rightarrow\; f(\xi_1, u_1) - f(\xi_2, u_2) \in \mathcal{T}_{\xi_1 - \xi_2} K \,. \tag{3}$$

The proof will be given in Section 3.5; first, we discuss the applicability of this test, and we develop several technical results.

## 3.1 Orthants

We start by looking at a special case, namely $K = \mathbb{R}^n_{\geq 0}$ and $K_u = \mathbb{R}^m_{\geq 0}$ (with $\mathbb{B}_\mathcal{U} = \mathbb{R}^m$). Such systems are called *cooperative systems*.

The boundary points of $K$ are those points for which some coordinate is zero, so "$\xi_1 - \xi_2 \in \partial K$" means that $\xi_1 \succeq \xi_2$ and $\xi_1^i = \xi_2^i$ for at least one $i \in \{1, \ldots, n\}$. On the other hand, if $\xi_1 \succeq \xi_2$ and $\xi_1^i = \xi_2^i$ for $i \in I$ and $\xi_1^i > \xi_2^i$ for $i \in \{1, \ldots, n\} \setminus I$, the tangent cone $\mathcal{T}_{\xi_1 - \xi_2} K$ consists of all those vectors $v = (v_1, \ldots, v_n) \in \mathbb{R}^n$ such that $v_i \geq 0$ for $i \in I$ and $v_i$ is arbitrary in $\mathbb{R}$ otherwise. Therefore, Property (3) translates into the following statement:

$$\xi_1 \succeq \xi_2 \text{ and } \xi_1^i = \xi_2^i \text{ and } u_1 \succeq u_2 \;\Rightarrow\; f^i(\xi_1, u_1) \geq f^i(\xi_2, u_2) \tag{4}$$

holding for all $i = 1, 2, \ldots n$, all $u_1, u_2 \in \mathcal{U}$, and all $\xi_1, \xi_2 \in \mathcal{V}$ (where $f^i$ denotes the $i$th component of $f$). In particular, for systems with no inputs

$$\dot{x} = f(x), \tag{5}$$

one recovers the well-known characterization for cooperativity (cf. [24]):

$$\xi_1 \succeq \xi_2 \text{ and } \xi_1^i = \xi_2^i \;\Rightarrow\; f^i(\xi_1) \geq f^i(\xi_2) \tag{6}$$

must hold for all $i = 1, 2, \ldots n$ and all $\xi_1, \xi_2 \in \mathcal{V}$.

When $X$ is strictly star-shaped, and in particular if $X$ is convex, cf. Remark 3.1, one could equally well require condition (4) to hold for all $\xi_1, \xi_2 \in X$. Indeed, pick any $\xi_1 \succeq \xi_2$, and suppose that $\xi_1^i = \xi_2^i$ for $i \in I$ and $\xi_1^i > \xi_2^i$ for $i \in \{1, \ldots, n\} \setminus I$. Pick sequences $\xi_1^k \to \xi_1$ and $\xi_2^k \to \xi_2$ so that, for all $k$, $\xi_1^k, \xi_2^k \in \mathcal{V}$, $\xi_1^k \succeq \xi_2^k$ and $(\xi_1^k)^i = (\xi_2^k)^i$ for $i \in I$ (this can be done by choosing an appropriate projection $L$ in Remark 3.1). Since the property holds for elements in $\mathcal{V}$, we have that $f^i(\xi_1^k, u_1) \geq f^i(\xi_2^k, u_2)$ for all $k = 1, 2, \ldots$ and all $i \in I$. By continuity, taking limits as $k \to \infty$, we also have then that $f^i(\xi_1, u_1) \geq f^i(\xi_2, u_2)$. On the other hand, if $\mathcal{U}$ also satisfies an approximability property, then by continuity one proves similarly that it is enough to check the condition (4) for $u_1, u_2$ belonging to the interior $\mathcal{W} = \text{int}\,\mathcal{U}$. In summary, we can say that if $X$ and $\mathcal{U}$ are both convex, then it is equivalent to check condition (4) for elements in the sets or in their respective interiors.

One can also rephrase the inequalities in terms of the partial derivatives of the components of $f$. Let us call a subset $S$ of an ordered Banach space *order-convex* ("p-convex" in [24]) if, for every $x$ and $y$ in $S$ with $x \succeq y$ and every $0 \leq \lambda \leq 1$, the element $\lambda x + (1-\lambda) y$ is in $S$. For instance, any convex set is order-convex, for all possible orders. We have the following easy fact, which generalizes Remark 4.1.1 in [24]:

**Proposition 3.3** Suppose that $\mathbb{B}_\mathcal{U} = \mathbb{R}^m$, $\mathcal{U}$ satisfies an approximability property, and both $\mathcal{V}$ and $\mathcal{W} = \text{int}\,\mathcal{U}$ are order-convex (for instance, these properties hold if both $\mathcal{V}$ and $\mathcal{U}$ are convex). Assume that $f$ is continuously differentiable. Then, the system (1) is cooperative if and only if the following properties hold:

$$\frac{\partial f^i}{\partial x^j}(x, u) \geq 0 \quad \forall\, x \in \mathcal{V},\; \forall\, u \in \mathcal{W},\; \forall\, i \neq j \tag{7}$$

and
$$\frac{\partial f^i}{\partial u^j}(x,u) \geq 0 \quad \forall\, x \in X,\ \forall\, u \in \mathcal{W} \tag{8}$$

for all $i \in \{1,2,\ldots n\}$ and all $j \in \{1,2,\ldots m\}$.

*Proof.* We will prove that these two conditions are equivalent to condition (4) holding for all $i = 1,2,\ldots n$, all $u_1, u_2 \in \mathcal{W}$, and all $\xi_1, \xi_2 \in \mathcal{V}$. Necessity does not require the order-convexity assumption. Pick any $\xi \in \mathcal{V}$, $u \in \mathcal{W}$, and pair $i \neq j$. We take $\xi_1 = \xi$, $u_1 = u_2 = u$, and $\xi_2(\lambda) = \xi + \lambda e_j$, where $e_j$ is the canonical basis vector having all coordinates $\neq j$ equal to zero and its $j$th coordinate one, with $\lambda < 0$ near enough to zero so that $\xi_2(\lambda) \in \mathcal{V}$. Notice that, for all such $\lambda$, $\xi_1 \succeq \xi_2(\lambda)$ and $\xi_1^i = \xi_2(\lambda)^i$ (in fact, $\xi_1^\ell = \xi_2(\lambda)^\ell$ for all $\ell \neq j$). Therefore condition (4) gives that $f^i(\xi_1, u) \geq f^i(\xi_2(\lambda), u)$ for all negative $\lambda \approx 0$. A similar argument shows that $f^i(\xi_1, u) \leq f^i(\xi_2(\lambda), u)$ for all positive $\lambda \approx 0$. Thus $f^i(\xi_2(\lambda), u)$ is increasing in a neighborhood of $\lambda = 0$, and this implies Property (7). A similar argument establishes Property (8).

For the converse, as in [24], we simply use the Fundamental Theorem of Calculus to write

$$f^i(\xi_2, u_1) - f^i(\xi_1, u_1) = \int_0^1 \sum_{j=1}^n \frac{\partial f^i}{\partial x^j}(\xi_1 + r(\xi_2 - \xi_1), u_1)\,(\xi_2^j - \xi_1^j)\, dr \tag{9}$$

and

$$f^i(\xi_2, u_2) - f^i(\xi_2, u_1) = \int_0^1 \sum_{j=1}^m \frac{\partial f^i}{\partial u^j}(\xi_2, u_2 + r(u_2^j - u_1^j))\,(u_2^j - u_1^j)\, dr \tag{10}$$

for any $i = 1,2,\ldots n$, $u_1, u_2 \in \mathcal{W}$, and $\xi_1, \xi_2 \in \mathcal{V}$. Pick any $i \in \{1,2,\ldots n\}$, $u_1, u_2 \in \mathcal{W}$, and $\xi_1, \xi_2 \in \mathcal{V}$, and suppose that $\xi_1 \succeq \xi_2$, $\xi_1^i = \xi_2^i$, and $u_1 \succeq u_2$. We need to show that $f^i(\xi_2, u_2) \leq f^i(\xi_1, u_1)$. Since the integrand in (9) vanishes when $j = i$, and also $\partial f^i/\partial x^j \geq 0$ and $\xi_2^j - \xi_1^j \leq 0$ for $j \neq i$, it follows that $f^i(\xi_2, u_1) \leq f^i(\xi_1, u_1)$. Similarly, Formula (10) gives us that $f^i(\xi_2, u_2) \leq f^i(\xi_2, u_1)$, completing the proof. ∎

For systems (5) without inputs, Property (7) is the well-known characterization "$\frac{\partial f^i}{\partial x^j} \geq 0$ for all $i \neq j$" of cooperativity. Interestingly, the authors of [22] use this property, for systems as in (1) but where inputs $u$ are seen as constant parameters, as a definition of (parameterized) cooperative systems, but monotonicity with respect to time-varying inputs is not exploited there. The terminology "cooperative" is motivated by this property: the different variables $x^i$ have a positive influence on each other.

### Other orthants

More general orthants can be treated by the trick used in Section 3.5 in [24]. Any orthant $K$ in $\mathbb{R}^n$ has the form
$$K^{(\varepsilon)} = \{x \in \mathbb{R}^n \mid (-1)^{\varepsilon_i} x_i \geq 0,\ i = 1,\ldots, n\}$$
for some binary vector $\varepsilon = (\varepsilon_1, \ldots, \varepsilon_n) \in \{0,1\}^n$. Note that $K^{(\varepsilon)} = P\mathbb{R}_{\geq 0}^n$, where $P : \mathbb{R}^n \to \mathbb{R}^n$ is the linear mapping given by the matrix $P = \mathrm{diag}\,((-1)^{\varepsilon_1}, \ldots, (-1)^{\varepsilon_n})$. Similarly, if the cone $K_u$ defining the order for $\mathcal{U}$ is an orthant $K^{(\delta)}$, we can view it as $Q\mathbb{R}_{\geq 0}^m$, for a similar map $Q = \mathrm{diag}\,((-1)^{\delta_1}, \ldots, (-1)^{\delta_m})$. Monotonicity of $\dot{x} = f(x,u)$ under these orders is equivalent to monotonicity of $\dot{z} = g(z,v)$, where $g(z,v) = Pf(Pz, Qv)$, under the already studied orders given by $\mathbb{R}_{\geq 0}^n$ and $\mathbb{R}_{\geq 0}^m$. This is because the change of variables $z(t) = Px(t)$, $v(t) = Qu(t)$ transforms solutions of one system into the other (and viceversa), and both $P$ and $Q$ preserve the respective orders ($\xi_1 \succeq \xi_2$ is equivalent to $(P\xi_1)^i \geq (P\xi_2)^i$ for all $i \in \{1, \ldots, n\}$, and similarly for input values). Thus we conclude:

**Corollary 3.4** Under the assumptions in Proposition 3.3, and for the orders induced from orthants $K^{(\varepsilon)}$ and $K^{(\delta)}$, the system (1) is monotone if and only if the following properties hold:

$$(-1)^{\varepsilon_i + \varepsilon_j} \frac{\partial f^i}{\partial x^j}(x, u) \geq 0 \quad \forall\, x \in \mathcal{V},\ \forall\, u \in \mathcal{W},\ \forall\, i \neq j \tag{11}$$

and
$$(-1)^{\varepsilon_i+\delta_j} \frac{\partial f^i}{\partial u^j}(x,u) \geq 0 \quad \forall\, x \in X,\ \forall\, u \in \mathcal{W} \tag{12}$$
for all $i \in \{1, 2, \ldots n\}$ and all $j \in \{1, 2, \ldots m\}$. □

## 3.2 An example

Let us clarify the above definitions and notations with an example. We consider the partial order $\succeq$ obtained by letting $K = \mathbb{R}_{\leq 0} \times \mathbb{R}_{\geq 0}$. Using the previous notations, we can write this as $K = K^{(\varepsilon)}$, where $\varepsilon = (1, 0)$. We will consider the input space $\mathcal{U} = \mathbb{R}_{\geq 0}$, with the standard ordering in $\mathbb{R}$ (i.e., $K_u = \mathbb{R}_{\geq 0}$, or $K_u = K^{(\delta)}$ with $\delta = (0)$).

Observe that the boundary points of this cone $K$ are those points of the forms $p = (0, a)$ or $q = (-a, 0)$, for some $a \geq 0$, and the tangent cones are respectively $\mathcal{T}_p K = \mathbb{R}_{\leq 0} \times \mathbb{R}$ and $\mathcal{T}_p K = \mathbb{R} \times \mathbb{R}_{\geq 0}$, see Fig. 1. Under the assumptions of Corollary 3.4, a system is monotone with respect to these orders

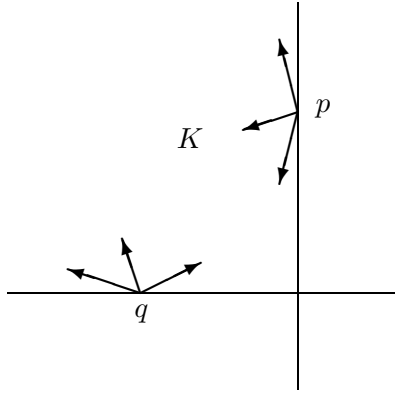

Figure 1: Example of cone and tangents

if and only if the following four inequalities hold everywhere:
$$\frac{\partial f^1}{\partial x^2} \leq 0,\ \frac{\partial f^2}{\partial x^1} \leq 0,\ \frac{\partial f^1}{\partial u} \leq 0,\ \frac{\partial f^2}{\partial u} \geq 0\,.$$

A special class of systems of this type is afforded by systems as follows:
$$\left. \begin{array}{rcl} \dot{x}^1 &=& -u\theta_1(x^1) + \theta_2(1 - x^1 - x^2) \\ \dot{x}^2 &=& u\theta_3(1 - x^1 - x^2) - \theta_4(x^2) \end{array} \right\} := f(x^1, x^2, u) \tag{13}$$

where the functions $\theta_i$'s are $\mathcal{C}^1$ increasing functions such that $\theta_i(0) = 0$. The system is regarded as evolving on the triangle
$$X = \Delta := \{[x, y] : x \geq 0,\ y \geq 0, x + y \leq 1\}$$

which is easily seen to be invariant for the dynamics. The following fact is immediate from the above discussion:

**Lemma 3.5** *The system (13) is monotone with respect to the given orders.* □

## 3.3 Competitive systems

A related class of dynamical systems are *competitive systems*; these are usually defined, for differential equations as in (5), by the condition $\frac{\partial f_j}{\partial x^i} \leq 0$ for all $i \neq j$. One way of understanding such a condition is by noticing that a system is competitive if and only if its "reversed time" version, $\dot{x} = -f(x)$, is cooperative. More generally, for generic partial orders $\succeq$, we can think of systems with "competitive dynamics" w.r.t $\succeq$, according to the following definition:

$$x_1 \succeq x_2 \Rightarrow \phi(t, x_1) \succeq \phi(t, x_2) \qquad \forall t \leq 0. \tag{14}$$

We point out that systems with "competitive dynamics" need not be monotone dynamical systems (in forward time); therefore, all the "good" properties concerning $\omega$-limit sets and convergence to attractors mentioned earlier are not necessarily true (although analogous statements are possible concerning their $\alpha$-limit sets and asymptotic behavior of trajectories as $t \to -\infty$).

As it turns out to be useful in applications, we explicitly give the definition of *controlled competitive system*:

**Definition 3.6** A controlled dynamical system $\phi : \mathbb{R}_{\geq 0} \times X \times \mathcal{C}_\mathcal{U} \to X$ is *monotone in backward time* (or has *competitive dynamics*) if the implication below holds:

$$u_1 \succeq u_2, \ x_1 \succeq x_2 \Rightarrow \phi(t, x_1, u_1) \succeq \phi(t, x_2, u_2) \qquad \forall t \leq 0. \qquad \square$$

The following infinitesimal characterization of "competitivity" follows looking at (3) backwards in time:

$$\xi_1 - \xi_2 \in \partial K \text{ and } u_1 \succeq u_2 \ \Rightarrow \ f(\xi_2, u_2) - f(\xi_1, u_1) \in \mathcal{T}_{\xi_1 - \xi_2} K$$

In particular, when $K$ and $K_u$ are the $n$ and $m$ dimensional positive orthants, and under the assumptions in Proposition 3.3, we obtain the following differential characterization of controlled competitive systems:

$$\frac{\partial f^i}{\partial x^j}(x, u) \leq 0 \quad \forall x \in X, \ \forall u \in \mathcal{U} \quad \forall i \neq j$$

together with

$$\frac{\partial f^i}{\partial u^j}(x, u) \leq 0 \quad \forall x \in X, \ \forall u \in \mathcal{U}$$

for all $i \in \{1, 2, \ldots n\}$ and all $j \in \{1, 2, \ldots m\}$.

### 3.4 Preliminaries on invariance

We now return to the proof of Theorem 1.

**Lemma 3.7** The set $\mathcal{V}$ is forward invariant for (1), i.e., for each $\xi \in \mathcal{V}$ and each $u \in \mathcal{U}_\infty$, $\phi(t, \xi, u) \in \mathcal{V}$ for all $t \geq 0$.

*Proof.* Pick any $\xi \in \mathcal{V}$, $u \in \mathcal{U}_\infty$, and $t_0 \geq 0$. Viewing (1) as a system defined on an open set of states $\widetilde{X}$ which contains $X$, we consider the mapping $\alpha : \mathcal{V} \to \mathbb{B}$ given by $\alpha(x) = \phi(t_0, x, u)$ (with the same $u$ and $t_0$). The image of $\alpha$ must contain a neighborhood $W$ of $\xi' = \alpha(\xi)$; see e.g. Lemma 4.3.8 in [26]. Thus, $W \subseteq X$, which means that $\xi' \in \text{int } X$, as desired. ∎

**Remark 3.8** The converse of Lemma 3.7 is also true, namely, if $\dot{x} = f(x, u)$ is a system defined on some neighborhood $\widetilde{X}$ of $X$ and if $\mathcal{V} = \text{int } X$ is forward invariant under solutions of this system, then $X$ is itself invariant. To see this, pick any $\xi \in X$ and a sequence $\xi^i \to \xi$ of elements of $\mathcal{V}$. For any $t$, $i$, and $u$, $\phi(t, \xi^i, u) \in \mathcal{V}$, so $\phi(t, \xi, u) = \lim_{i \to \infty} \phi(t, \xi^i, u) \in \text{clos } \mathcal{V} = X$. $\square$

We introduce the following subset of $X \times X$:

$$\Gamma := \{(\xi_1, \xi_2) \in X \times X \mid \xi_1 \succeq \xi_2\}$$

(we use row-vector notation, for notational simplicity) as well as the set

$$\mathcal{U}^{[2]} := \{(u_1, u_2) \in \mathcal{U} \times \mathcal{U} \mid u_1 \succeq u_2\}$$

and observe that $\Gamma$ is a closed subset since the defining condition means that $\xi_1 - \xi_2 \in K$, and $K$ is closed, and similarly $\mathcal{U}^{[2]}$ is also closed. Next we consider the system with state-space $X \times X$ and input-value set $\mathcal{U}^{[2]}$ whose dynamics

$$\dot{x} = f^{[2]}(x, u) \tag{15}$$

are given, in block form using $x = (x_1, x_2) \in X \times X$ and $u = (u_1, u_2) \in \mathcal{U}^{[2]}$, as:

$$\begin{aligned} \dot{x}_1 &= f(x_1, u_1) \\ \dot{x}_2 &= f(x_2, u_2) \end{aligned}$$

(two copies of the same system, driven by the different $u_i$'s). We denote by $\mathcal{U}^{[2]}_\infty$ the set of all possible inputs to this composite system, i.e. the set of all Lebesgue-measurable locally essentially bounded functions $u: [0, \infty) \to \mathcal{U}^{[2]}$. Since by Lemma 3.7 the interior $\mathcal{V}$ of $X$ is forward invariant for (1), it holds that $\phi^{[2]}(t, \xi, u)$ belongs to $\mathcal{V} \times \mathcal{V}$ whenever $\xi \in \mathcal{V} \times \mathcal{V}$ and $u \in \mathcal{U}^{[2]}_\infty$.

Observe that the definition of monotonicity amounts to the requirement that: *for each $\xi \in \Gamma$, and each $u \in \mathcal{U}^{[2]}_\infty$, the solution $\phi^{[2]}(t, \xi, u)$ of (15) with initial condition $x(0) = \xi$ belongs to $\Gamma$ for all $t \geq 0$* (forward invariance of $\Gamma$ with respect to (15)).

The set:
$$\Gamma_0 := \Gamma \bigcap (\mathcal{V} \times \mathcal{V})$$

is closed relative to $\mathcal{V} \times \mathcal{V}$. The following elementary remark will be very useful:

**Lemma 3.9** The system (1) is monotone if and only if the set $\Gamma_0$ is forward invariant for the system (15) restricted to $\mathcal{V} \times \mathcal{V}$.

*Proof.* We must show that monotonicity is the same as:

$$\xi \in \Gamma_0 \text{ and } u \in \mathcal{U}^{[2]}_\infty \quad \Rightarrow \quad \phi^{[2]}(t, \xi, u) \in \Gamma_0 \ \forall t \geq 0. \tag{16}$$

Necessity is clear, since if the system is monotone then $\phi^{[2]}(t, \xi, u) \in \Gamma$ holds for all $\xi \in \Gamma \supseteq \Gamma_0$ and all $t \geq 0$, and we already remarked that $\phi^{[2]}(t, \xi, u) \in \mathcal{V} \times \mathcal{V}$ whenever $\xi \in \Gamma_0$. Conversely, suppose that (16) holds. Pick any $\xi \in \Gamma$. The approximability hypothesis provides a sequence $\{\xi^i\} \subseteq \Gamma_0$ such that $\xi^i \to \xi$ as $i \to \infty$. Fix any $u \in \mathcal{U}^{[2]}_\infty$ and any $t \geq 0$. Then $\phi^{[2]}(t, \xi^i, u) \in \Gamma_0 \subseteq \Gamma$ for all $i$, so taking limits and using continuity of $\phi^{[2]}$ on initial conditions gives that $\phi^{[2]}(t, \xi, u) \in \Gamma$, as required. ∎

**Lemma 3.10** For any $\xi = (\xi_1, \xi_2) \in \Gamma_0$ and any $u = (u_1, u_2) \in \mathcal{U}^{[2]}$, the following two properties are equivalent:

$$f(\xi_1, u_1) - f(\xi_2, u_2) \in \mathcal{T}_{\xi_1 - \xi_2} K \tag{17}$$

$$f^{[2]}(\xi, u) \in \mathcal{T}_\xi \Gamma_0. \tag{18}$$

*Proof.* Suppose that (17) holds, so there are sequences $t_i \searrow 0$ and $\{\eta^i\} \subseteq K$ such that $\eta^i \to \xi_1 - \xi_2$ and

$$\frac{1}{t_i} \left( \eta^i - (\xi_1 - \xi_2) \right) \to f(\xi_1, u_1) - f(\xi_2, u_2) \quad \text{as } i \to \infty. \tag{19}$$

Since $\mathcal{V}$ is open, the solution $x(t) = \phi(t, \xi_1, \bar{u})$ of $\dot{x} = f(x, \bar{u})$ with input $\bar{u} \equiv u_1$ and initial condition $x(0) = \xi_1$ takes values in $\mathcal{V}$ for all sufficiently small $t$. Thus, restricting to a subsequence, we may without loss of generality assume that $\xi^i_1 := x(t_i)$ is in $\mathcal{V}$ for all $i$. Note that, by definition of solution,

$$\frac{1}{t_i} (\xi^i_1 - \xi_1) \to f(\xi_1, u_1) \quad \text{as } i \to \infty. \tag{20}$$

Subtracting (19) from (20), we obtain:

$$\frac{1}{t_i} (\xi^i_2 - \xi_2) \to f(\xi_2, u_2) \quad \text{as } i \to \infty \tag{21}$$

with $\xi^i_2 := \xi^i_1 - \eta^i$. Since $\xi^i_1 \to \xi_1$ and $\eta^i \to \xi_1 - \xi_2$ as $i \to \infty$, the sequence $\xi^i_2$ converges to $\xi_2 \in \mathcal{V}$. Using once again that $\mathcal{V}$ is open, we may assume without loss of generality that $\xi^i_2 \in \mathcal{V}$ for all $i$.

Moreover, $\xi_1^i - \xi_2^i = \eta^i \in K$, i.e., $\xi^i := (\xi_1^i, \xi_2^i) \in \Gamma$ for all $i$, which means that $\xi^i$ is in $\Gamma_0$ for all $i$, and, from (20) and (21):

$$\frac{1}{t_i}\left(\xi^i - \xi\right) \to (f(\xi_1, u_1), f(\xi_2, u_2)) \quad \text{as } i \to \infty \tag{22}$$

so that Property (18) is verified. Conversely, if Property (18) holds, then there are sequences $t_i \searrow 0$ and $\xi^i := (\xi_1^i, \xi_2^i) \in \Gamma_0$ such that (22) holds, and then, coordinatewise, we have both (20) and (21), which subtracted and defining $\eta^i := \xi_1^i - \xi_2^i$ give (19); this establishes Property (17). ∎

### 3.5 Proof of Theorem 1

Suppose that the system (1) is monotone, and fix any input-value pair $u^0 = (u_1^0, u_2^0) \in \mathcal{U}^{[2]}$. Lemma 3.9 says that the set $\Gamma_0$ is forward invariant for the system (15) restricted to $\mathcal{V} \times \mathcal{V}$. This implies, in particular, that every solution of the differential equation

$$\dot{x} = f^{[2]}(x, u^0)$$

with $x(0) \in \Gamma_0$ remains in $\Gamma_0$ for all $t \geq 0$ (where we think of $u^0$ as a constant input). We may view this differential equation as a (single-valued) differential inclusion $\dot{x} \in F(x)$ on $\mathcal{V} \times \mathcal{V}$, where $F(\xi) = \{f^{[2]}(\xi, u^0)\}$, for which the set $\Gamma_0$ is strongly invariant. Thus, Theorem 3 in the Appendix implies that $F(\xi) \subseteq \mathcal{T}_\xi \Gamma_0$ for all $\xi \in \Gamma_0$. In other words, Property (18) holds, at all $\xi \in \Gamma_0$, for the given $u = u^0$, and hence, by Lemma 3.10, $f(\xi_1, u_1^0) - f(\xi_2, u_2^0) \in \mathcal{T}_{\xi_1 - \xi_2} K$ for all $(\xi_1, \xi_2) \in \Gamma_0$ and this $u^0$. As $u^0$ was an arbitrary element of $\mathcal{U}^{[2]}$, Property (2) follows.

Conversely, suppose that (2) holds. By Lemma 3.10, we know that Property (18) holds for all $(\xi_1, \xi_2) \in \Gamma_0$ and all $(u_1, u_2) \in \mathcal{U}^{[2]}$. To show monotonicity of the system (1), we need to prove that $\Gamma_0$ is invariant for the system (15) when restricted to $\mathcal{V} \times \mathcal{V}$. So pick any $\xi^0 \in \Gamma_0$, any $u^0 \in \mathcal{U}_\infty^{[2]}$, and any $t^0 > 0$; we must prove that $\phi^{[2]}(t^0, \xi^0, u^0) \in \Gamma_0$. The input function $u^0$ being locally bounded means that there is some compact subset $C \subseteq \mathcal{U}$ such that $u(t)$ belongs to the compact subset $\mathcal{U}_C^{[2]} = \mathcal{U}^{[2]} \bigcap C \times C$ of $\mathbb{B}_\mathcal{U} \times \mathbb{B}_\mathcal{U}$, for (almost) all $t \in [0, t^0]$. We introduce the following compact-valued, locally bounded, and locally Lipschitz set-valued function:

$$F_C(\xi) := \left\{f^{[2]}(\xi, u) \mid u \in \mathcal{U}_C^{[2]}\right\}$$

on $\mathcal{V} \times \mathcal{V}$. We already remarked that Property (17) holds, i.e., $\{f^{[2]}(\xi, u) \mid u \in \mathcal{U}^{[2]}\} \subseteq \mathcal{T}_\xi \Gamma_0$, for all $(\xi_1, \xi_2) \in \Gamma_0$, so it is true in particular that $F_C(\xi) \subseteq \mathcal{T}_\xi \Gamma_0$. Thus, Theorem 3 in the Appendix implies that $\Gamma_0$ is strongly invariant with respect to $F_C$. Thus, since $x(\cdot) = \phi^{[2]}(\cdot, \xi^0, u^0)$ restricted to $[0, t^0]$ satisfies $\dot{x} \in F_C(x)$, we conclude that $x(t^0) \in \Gamma_0$, as required.

Finally, we show that (2) and (3) are equivalent. Since (3) is a particular case of (2), we only need to verify that $f(\xi_1, u_1) - f(\xi_2, u_2) \in \mathcal{T}_{\xi_1 - \xi_2} K$ when $\xi_1 - \xi_2 \in \text{int } K$. This is a consequence of the general fact that $\mathcal{T}_\xi S = \mathbb{R}^n$ whenever $\xi$ is in the interior of a set $S$. ∎

## 4 Cascades of monotone systems

Cascade structures with triangular form

$$\begin{aligned}
\dot{x}_1 &= f_1(x_1, x_2, \ldots, x_N, u) \\
\dot{x}_2 &= f_2(x_2, \ldots, x_N, u) \\
&\vdots \quad \vdots \\
\dot{x}_i &= f_i(x_i, \ldots, x_N, u) \\
&\vdots \quad \vdots \\
\dot{x}_N &= f_N(x_N, u)
\end{aligned} \tag{23}$$

are of special interest. A simple sufficient condition for monotonicity of systems (23) is as follows.

**Proposition 4.1** Assume that there exist positivity cones $K_1, K_2, \ldots, K_{N+1}$ (of suitable dimensions) so that each of the $x_i$-subsystems in (23) is a controlled monotone dynamical system with respect to the $K_i$-induced partial order (as far as states are concerned) and with respect to the $K_{i+1}, \ldots, K_{N+1}$-induced partial orders as far as inputs are concerned. Then, the overall cascaded interconnection (23) is monotone with respect to the order induced by the positivity cone $K_1 \times K_2 \times \ldots \times K_N$ on states and $K_{N+1}$ on inputs.

*Proof.* We first prove the result for the case $N = 2$:

$$\begin{aligned} \dot{x}_1 &= f_1(x_1, x_2, u) \\ \dot{x}_2 &= f_2(x_2, u) \,. \end{aligned} \quad (24)$$

Let $\succeq_1$ and $\succeq_2$ be the partial orders induced by the cones $K_1, K_2$ and $\succeq_u$ on inputs. Pick any two inputs $u^a \succeq_u u^b$. By hypothesis we have, for each two states $\xi^a = (\xi_1^a, \xi_2^a)$ and $\xi^b = (\xi_1^b, \xi_2^b)$:

$$\xi_2^a \succeq_2 \xi_2^b \Rightarrow \phi_2(t, \xi_2^a, u^a) \succeq_2 \phi_2(t, \xi_2^b, u^b) \quad \forall\, t \geq 0 \quad (25)$$

as well as, for all functions $x_2^a, x_2^b$:

$$\xi_1^a \succeq_1 \xi_1^b, \ x_2^a \succeq_2 x_2^b \Rightarrow \phi_1(t, \xi_1^a, x_2^a, u^a) \succeq_1 \phi_1(t, \xi_1^b, x_2^b, u^b) \quad \forall\, t \geq 0. \quad (26)$$

Combining (25) and (26) yields

$$\begin{aligned} \xi_a^1 \succeq_1 \xi_1^b \ \& \ \xi_2^a \succeq_2 \xi_2^b \ &\Rightarrow\ \xi_a^1 \succeq_1 \xi_1^b \ \& \ \phi_2(t, \xi_2^a, u^a) \succeq_2 \phi_2(t, \xi_2^b, u^b) \quad \forall\, t \geq 0 \\ &\Rightarrow\ \phi_1(t, \xi_1^a, x_2^a, u^a) \succeq_1 \phi_1(t, \xi_1^b, x_2^b, u^b) \ \& \ \phi_2(t, \xi_2^a, u^a) \succeq_2 \phi_2(t, \xi_2^b, u^b) \quad \forall\, t \geq 0 \,. \end{aligned} \quad (27)$$

Then, defining $K := K_1 \times K_2$ and letting $\succeq$ denote the corresponding partial order, implication (27) reads:

$$\begin{bmatrix} \xi_1^a \\ \xi_2^a \end{bmatrix} \succeq \begin{bmatrix} \xi_1^b \\ \xi_2^b \end{bmatrix} \Rightarrow \phi\left(t, \begin{bmatrix} \xi_1^a \\ \xi_2^a \end{bmatrix}, u^a\right) \succeq \phi\left(t, \begin{bmatrix} \xi_1^b \\ \xi_2^b \end{bmatrix}, u^b\right) \quad \forall\, t \geq 0. \quad (28)$$

The proof for arbitrary $N$ follows by induction. ∎

## 5 Static Input/State and Input/Output Characteristics

A notion of "Cauchy gain" was introduced in [27] to quantify amplification of signals in a manner useful for biological applications. For monotone dynamical systems satisfying an additional property, it is possible to obtain tight estimates of Cauchy gains. This is achieved by showing that the output values $y(t)$ corresponding to an input $u(\cdot)$ are always "sandwiched" in between the outputs corresponding to two constant inputs which bound the range of $u(\cdot)$. This additional property motivated our looking at monotone systems to start with; we now start discussion of that topic.

**Definition 5.1** We say that a controlled dynamical system (1) is endowed with the *static Input/State characteristic*

$$k_x(\cdot) : \mathcal{U} \to X$$

if for each constant input $u(t) \equiv \bar{u}$ there exists a (necessarily unique) globally asymptotically stable equilibrium $k_x(\bar{u})$. For systems with an output map $y = h(x)$, we also define the *static Input/Output characteristic* as $k_y(\bar{u}) := h(k_x(\bar{u}))$, provided that an Input/State characteristic exists and that $h$ is continuous. □

The paper [22] (see also citemuratori-rinaldi for linear systems) provides very useful results which can be used to show the existence of I/S characteristics, for cooperative systems with scalar inputs and whose state space is the positive orthant, and in particular to the study of the question of when $k_x(\bar{u})$ is strictly positive.

**Remark 5.2** Observe that, if the system (1) is monotone and it admits a static Input/State characteristic $k_x$, then $k_x$ must be nondecreasing with respect to the orders in question: $\bar{u} \succeq \bar{v}$ in $\mathcal{U}$ implies $k_x(\bar{u}) \succeq k_x(\bar{v})$. Indeed, given any initial state $\xi$, monotonicity says that $\phi(t, \xi, u) \succeq \phi(t, \xi, v)$ for all $t$, where $u(t) \equiv \bar{u}$ and $v(t) \equiv \bar{v}$. Taking limits as $t \to \infty$ gives the desired conclusion. □

**Remark 5.3** (Continuity of $k_x$) Suppose that for a system (1) there is a map $k_x : \mathcal{U} \to X$ with the property that $k_x(\bar{u})$ is the unique steady state of the system $\dot{x} = f(x, \bar{u})$ (constant input $u \equiv \bar{u}$). When $k_x(\bar{u})$ is a globally asymptotically stable state for $\dot{x} = f(x, \bar{u})$, as is the case for I/S characteristics, it follows that the function $k_x$ must be continuous, see Proposition 5.7 below. However, continuity is always true provided only that $k_x$ be locally bounded, i.e. that $k_x(V)$ is a bounded set whenever $V \subseteq \mathcal{U}$ is compact. This is because $k_x$ has a closed graph, since $k_x(\bar{u}) = \bar{x}$ means that $f(\bar{x}, \bar{u}) = 0$, and any locally bounded map with a closed graph (in finite-dimensional spaces) must be continuous. (Proof: suppose that $\bar{u}_i \to \bar{u}$, and consider the sequence $\bar{x}_i = k_x(\bar{u}_i)$; by local boundedness, it is only necessary to prove that every limit point of this sequence equals $k_x(\bar{u})$. So suppose that $\bar{x}_{i_j} \to \bar{x}'$; then $(\bar{u}_{i_j}, \bar{x}_{i_j}) \to (\bar{u}, \bar{x}')$, so by the closedness of the graph of $k_x$ we know that $(\bar{u}, \bar{x}')$ belongs to its graph, and thus $\bar{x}' = \bar{x}$, as desired.) Therefore, local boundedness, and hence continuity of $k_x$, would follow if one knows that $k_x$ is monotone, so that $k([a, b])$ is always bounded, even if the stability condition does not hold, at least if the order is "reasonable" enough, as in the next definition. Note that $k_y$ is continuous whenever $k_x$ is, since the output map $h$ has been assumed to be continuous. □

Under weak assumptions, existence of a static Input/State characteristic implies that the system behaves well with respect to arbitrary bounded inputs as well as inputs that converge to some limit. For convenience in stating results along those lines, we introduce the following terminology.

**Definition 5.4** The order on $X$ is *bounded* if the following two properties hold:

1. For each bounded subset $S \subseteq X$, there exist two elements $a, b \in \mathbb{B}$ such that $S \subseteq [a, b] = \{x \in X : a \preceq x \preceq b\}$.
2. For each $a, b \in \mathbb{B}$, the set $[a, b]$ is bounded. □

**Remark 5.5** Boundedness is a very mild assumption. In general, Property 1 holds if (and only if) $K$ has a nonempty interior, and Property 2 is a consequence of $K \cap -K = \{0\}$. To prove the first claim, suppose given a bounded set $S \subseteq X$, and take any element $\xi_0 \in \operatorname{int} K$. Now pick any $\varepsilon > 0$ such that the function $\alpha : \mathbb{B} \to \mathbb{B} : \xi \mapsto \varepsilon \xi + \xi_0$ maps $S \bigcup -S$ into the interior of $K$ (e.g., pick $\varepsilon = \varepsilon_0/c$, if the $\varepsilon_0$-ball around $\xi_0$ is included in $K$ and if each element of $S \bigcup -S$ has norm at most $c$). Then $\varepsilon \xi \succeq -\xi_0$ and $\xi_0 \succeq \varepsilon \xi$ for all $\xi \in S$, which means that $S \subseteq [-\xi_0/\varepsilon, \xi_0/\varepsilon]$. Conversely, suppose that Property 1 holds. Since $X$ has nonempty interior, there exist $\varepsilon > 0$ and $\xi_0 \in X$ so that $S := \varepsilon B + \xi_0 \subseteq X$, where $B$ is the open unit ball in $\mathbb{B}$. Suppose that there would be some $b$ such that $b \succeq \xi$ for all $\xi \in S$. This implies that $\varepsilon B + \eta \subseteq K$, where $\eta := b - \xi_0$. Since $\varepsilon B + \eta$ is an open set, this means that $K$ has a nonempty interior.

Next we show that the property $K \cap -K = \{0\}$ implies that $[a, b]$ is a bounded set, for any $a, b$. Assume by contradiction that a sequence $x_m \in [a, b]$ exists such that $|x_m| \to +\infty$ as $m \to +\infty$. Note that $x_m - a$ and $b - x_m$ both belong to $K$ for all $m$. Using compactness of the unit sphere, it is possible to take a subsequence $x_{m_k}$ so that the sequences $b_k := (b - x_{m_k})/|b - x_{m_k}|$ and $a_k := (x_{m_k} - a)/|a - x_{m_k}|$ both have limits as $k \to +\infty$. Since $K$ is a cone, $b_k \in K$ and $a_k \in K$ for each $k$, so that closedness of $K$ yields $\lim_k b_k := \bar{b} \in K$ and $\lim_k a_k := \bar{a} \in K$. Moreover, by continuity of norms $|\cdot|$, $|\bar{b}| = |\bar{a}| = 1$. A simple calculation exploiting the fact that $|x_{m_k}| \to +\infty$ allows us to conclude that $\bar{b} + \bar{a} = 0$. Then, the property $K \cap -K = \{0\}$ implies $\bar{a} = \bar{b} = 0$, contradicting the fact that $|\bar{a}| = 1$. □

**Proposition 5.6** Consider a monotone system (1) which is endowed with a static Input/State characteristic, and suppose that the order on the state space $X$ is bounded. Pick any input $u$ all whose values $u(t)$ lie in some interval $[c, d] \subseteq \mathcal{U}$. (For example, $u$ could be any bounded input, if $K$ is an orthant in $\mathbb{R}^n$, or more generally if the order in $\mathcal{U}$ is bounded.) Let $x(t) = \phi(t, \xi, u)$ be any trajectory of the system corresponding to this control. Then $\{x(t), t \geq 0\}$ is a bounded subset of $X$.

*Proof.* Let $x_1(t) = \phi(t, \xi, d)$, so $x_1(t) \to k_x(d)$ as $t \to \infty$ and, in particular, $x_1(\cdot)$ is bounded; so (bounded order), there is some $b \in \mathbb{B}$ such that $x_1(t) \preceq b$ for all $t \geq 0$. By monotonicity,

$$x(t) = \phi(t, \xi, u) \preceq \phi(t, \xi, d) = x_1(t) \preceq b$$

for all $t \geq 0$. A similar argument using the lower bound $c$ on $u$ shows that there is some $a \in \mathbb{B}$ such that $a \preceq x(t)$ for all $t$. Thus $x(t) \in [a, b]$ for all $t$, which implies, again appealing to the bounded order hypothesis, that $x(\cdot)$ is bounded. ∎

Certain standard facts concerning the robustness of stability will be useful. We collect the necessary results in the next statements, for easy reference.

**Proposition 5.7** If (1) is a monotone system which is endowed with a static Input/State characteristic $k_x$, then $k_x$ is a continuous map. Moreover for each $\bar{u} \in \mathcal{U}$, $\bar{x} = k_x(\bar{u})$, the following properties hold:

1. For each neighborhood $P$ of $\bar{x}$ in $X$ there exist a neighborhood $P_0$ of $\bar{x}$ in $X$, and a neighborhood $Q_0$ of $\bar{u}$ in $\mathcal{U}$, such that $\phi(t, \xi, u) \in P$ for all $t \geq 0$, all $\xi \in P_0$, and all inputs $u$ such that $u(t) \in Q_0$ for all $t \geq 0$.

2. If in addition the order on the state space $X$ is bounded, then, for each input $u$ all whose values $u(t)$ lie in some interval $[c, d] \subseteq \mathcal{U}$ and with the property that $u(t) \to \bar{u}$, and all initial states $\xi \in X$, necessarily $x(t) = \phi(t, \xi, u) \to \bar{x}$ as $t \to \infty$.

*Proof.* Consider any trajectory $x(t) = \phi(t, \xi, u)$ as in Property 2. By Proposition 5.6, we know that there is some compact $C \subseteq \mathbb{B}$ such that $x(t) \in C$ for all $t \geq 0$. Since $X$ is closed, we may assume that $C \subseteq X$. We are therefore in the following situation: the autonomous system $\dot{x} = f(x, \bar{u})$ admits $\bar{x}$ as a globally asymptotically stable equilibrium (with respect to the state space $X$) and the trajectory $x(\cdot)$ remains in a compact subset of the domain of attraction (of $\dot{x} = f(x, \bar{u})$ seen as a system on an open subset of $\mathbb{B}$ which contains $X$). The "converging input converging state" property then holds for this trajectory; see [25] or the "local ISS" results in [29] (the paper [25] was written for globally asymptotically stable systems in $\mathbb{R}^n$, and for $\mathcal{U}$ Euclidean, but these hypotheses were not used, and the result holds under the conditions that apply in this paper; see [28], Theorem 1, for details).

Property 1 is a consequence of the same results (once again, see [28], Theorem 1, for details).

The continuity of $k_x$ is a consequence of Property 1. As discussed in Remark 5.3, we only need to show that $k_x$ is locally bounded, for which it is enough to show that for each $\bar{u}$ there is some neighborhood $Q_0$ of $\bar{u}$ and some compact subset $P$ of $X$ such that $k_x(\mu) \in P$ for all $\mu \in Q_0$. Pick any $\bar{u}$, and any compact neighborhood $P$ of $\bar{x} = k_x(\bar{u})$. By Property 1, there exist a neighborhood $P_0$ of $\bar{x}$ in $X$, and a neighborhood $Q_0$ of $\bar{u}$ in $\mathcal{U}$, such that $\phi(t, \xi, u_\mu) \in P$ for all $t \geq 0$ whenever $\xi \in P_0$ and $u_\mu(t) \equiv \mu$ with $\mu \in Q_0$. In particular, this implies that $k_x(\mu) = \lim_{t \to \infty} \phi(t, \bar{x}, u_\mu) \in P$, which is as required. ∎

**Corollary 5.8** Suppose that the system $\dot{x} = f(x, u)$ with output $y = h(x)$ is monotone and has static Input/State and Input/Output characteristics $k_x$, $k_y$, and that the system $\dot{z} = g(z, y)$ (with input value space equal to the output value space of the first system) has a static Input/State characteristic $k_z$, it is monotone, and the order on its state space $Z$ is bounded. Assume that the order on outputs $y$ is such that whenever $y$ is bounded, there are $c, d$ such that $y(t) \in [c, d]$ for all $t \geq 0$ (for instance, $K$ might be an orthant). Then the cascade system

$$\begin{aligned} \dot{x} &= f(x, u), \quad y = h(x) \\ \dot{z} &= g(z, y) \end{aligned}$$

is a monotone system which admits the static Input/State characteristic $\widetilde{k}(\bar{u}) = (k_x(\bar{u}), k_z(k_y(\bar{u})))$.

*Proof.* Pick any $\bar{u}$. We must show that $\widetilde{k}(\bar{u})$ is a globally asymptotically stable equilibrium of the cascade. Pick any initial state $(\xi, \zeta)$ of the composite system, and let $x(t) = \phi_x(t, \xi, \bar{u})$ (input constantly equal to $\bar{u}$), $y(t) = h(x(t))$, and $z(t) = \phi_z(t, \zeta, y)$. Notice that $x(t) \to \bar{x}$ and $y(t) = h(x(t)) \to \bar{y} = k_y(\bar{u})$, so viewing $y$ as an input to the second system and using Property 2 in Proposition 5.7, we have that $z(t) \to \bar{z} = k_z(k_y(\bar{u}))$. This establishes attractivity. To show stability, pick any neighborhoods $P_x$ and $P_z$ of $\bar{x}$ and $\bar{z}$ respectively. By Property 1 in Proposition 5.7, there are neighborhoods $P_0$ and $Q_0$ such that $\zeta \in P_0$ and $y(t) \in Q_0$ for all $t \geq 0$ imply $\phi_z(t, \zeta, y) \in P_z$ for all $t \geq 0$. Consider $P_1 := P_x \cap h^{-1}(Q_0)$, which is a neighborhood of $\bar{x}$, and pick any neighborhood $P_2$ of $\bar{x}$ with the property that $\phi(t, \xi, \bar{u}) \in P_1$ for all $\xi \in P_2$ and all $t \geq 0$ (stability of the equilibrium $\bar{x}$). Then, for all $(\xi, \zeta) \in P_2 \times P_0$, $x(t) = \phi_x(t, \xi, \bar{u}) \in P_1$ (in particular, $x(t) \in P_x$ for all $t \geq 0$, so $y(t) = h(x(t)) \in Q_0$, and hence also $z(t) = \phi_z(t, \zeta, y) \in Q_z$ for all $t \geq 0$. ∎

In analogy to what usually done for autonomous dynamical systems, we define the $\Omega$-limit set of any function $\alpha : [0, \infty) \to A$, where $A$ is a topological space (we will apply this to state-space solutions and to outputs) as follows:

$$\Omega[\alpha] := \{a \in A \mid \exists\, t_k \to +\infty \text{ s.t. } \alpha(t_k) \to a \text{ as } k \to +\infty\}$$

(in general, this set may be empty). For inputs $u \in \mathcal{U}_\infty$, we also introduce these two sets:

$$\mathcal{L}_\leq[u] := \{\mu \in \mathcal{U} \mid \exists\, t_k \to +\infty \text{ and } \mu_k \to \mu \text{ as } k \to +\infty, \mu_k \in \mathcal{U}, \text{ such that } u(t) \succeq \mu_k\ \forall\, t \geq t_k\}$$

and

$$\mathcal{L}_\geq[u] := \{\mu \in \mathcal{U} \mid \exists\, t_k \to +\infty \text{ and } \mu_k \to \mu \text{ as } k \to +\infty, \mu_k \in \mathcal{U}, \text{ such that } \mu_k \succeq u(t)\ \forall\, t \geq t_k\}.$$

These notations are motivated by the following special case.

Suppose that we consider a *SISO* (single-input single-output) system, by which we mean a system for which $\mathbb{B}_\mathcal{U} = \mathbb{R}$ and $\mathbb{B}_\mathcal{Y} = \mathbb{R}$, taken with the usual orders. Given any scalar bounded input $u(\cdot)$, we denote $u_{\inf} := \liminf_{t \to +\infty} u(t)$ and $u_{\sup} := \limsup_{t \to +\infty} u(t)$. Then, $u_{\inf} \in \mathcal{L}_\leq[u]$ and $u_{\sup} \in \mathcal{L}_\geq[u]$, as follows by definition of lim inf and lim sup. Similarly, both $\liminf_{t \to +\infty} y(t)$ and $\limsup_{t \to +\infty} y(t)$ belong to $\Omega[y]$, for any output $y$.

**Proposition 5.9** Consider a monotone system (1), with static I/S and I/O characteristics $k_x$ and $k_y$. Then, for each initial condition $\xi$ and each input $u$, the solution $x(t) = \phi(t, \xi, u)$ and the corresponding output $y(t) = h(x(t))$ satisfy:

$$k_x\left(\mathcal{L}_\leq[u]\right) \preceq \Omega[x] \preceq k_x\left(\mathcal{L}_\geq[u]\right)$$

and

$$k_y\left(\mathcal{L}_\leq[u]\right) \preceq \Omega[y] \preceq k_y\left(\mathcal{L}_\geq[u]\right).$$

*Proof.* Pick any $\xi, u$, and the corresponding $x(\cdot)$ and $y(\cdot)$, and any element $\mu \in \mathcal{L}_\leq[u]$. Let $t_k \to +\infty$, $\mu_k \to \mu$, with all $\mu_k \in \mathcal{U}$, and $u(t) \succeq \mu_k$ for all $t \geq t_k$. By monotonicity of the system, for $t \geq t_k$ we have:

$$x(t, \xi, u) = x(t - t_k, x(t_k, \xi, u), u) \succeq x(t - t_k, x(t_k, \xi, u), \mu_k). \tag{29}$$

In particular, if $x(s_\ell) \to \zeta$ for some sequence $s_\ell \to \infty$, it follows that

$$\zeta \succeq \lim_{\ell \to \infty} x(s_\ell - t_k, x(t_k, \xi, u), \mu_k) = k_x(\mu_k).$$

Next, taking limits as $k \to \infty$, and using continuity of $k_x$, this proves that $\zeta \succeq k_x(\mu)$. This property holds for every elements $\zeta \in \Omega[x]$ and $\mu \in \mathcal{L}_\leq[u]$, so we have shown that $k_x\left(\mathcal{L}_\leq[u]\right) \preceq \Omega[x]$. The remaining inequalities are all proved in an entirely analogous fashion. ∎

In particular we have:

**Proposition 5.10** Consider a monotone SISO system (1), with static I/S and I/O characteristics $k_x(\cdot)$ and $k_y(\cdot)$. Then, the I/S and I/O characteristics are nondecreasing, and for each initial condition $\xi$ and each bounded input $u(\cdot)$, the following holds:

$$k_y(u_{\inf}) \leq \liminf_{t \to +\infty} y(t, \xi, u) \leq \limsup_{t \to +\infty} y(t, \xi, u) \leq k_y(u_{\sup}). \tag{30}$$

If, instead, outputs are ordered by $\geq$, then the I/O static characteristic is nonincreasing, and for each initial condition $\xi$ and each bounded input $u(\cdot)$, the following inequality holds:

$$k_y(u_{\sup}) \leq \liminf_{t \to +\infty} y(t, \xi, u) \leq \limsup_{t \to +\infty} y(t, \xi, u) \leq k_y(u_{\inf}). \tag{31}$$

*Proof.* The proof of the first statement is immediate from Proposition 5.9 and the properties: $u_{\inf} \in \mathcal{L}_{\leq}[u]$, $u_{\sup} \in \mathcal{L}_{\geq}[u]$, $\liminf_{t \to +\infty} y(t) \in \Omega[y]$, and $\limsup_{t \to +\infty} y(t) \in \Omega[y]$, and the second statement is proved in a similar fashion. ∎

**Remark 5.11** It is an immediate consequence of Proposition 5.10 that, if a monotone system admits a static I/O characteristic $k$, and if there is a class-$\mathcal{K}_\infty$ function $\gamma$ such that $|k(u) - k(v)| \leq \gamma(|u-v|)$ for all $u, v$ (for instance, if $k$ is Lipschitz with constant $\rho$ one may pick as $\gamma$ the linear function $\gamma(r) = \rho r$) then the system has a Cauchy gain (in the sense of [27]) $\gamma$ on bounded inputs. □

# 6 Feedback Interconnections of Monotone Systems

In this section, we study the stability of SISO monotone dynamical systems connected in feedback as in Fig. 2. Observe that such interconnections need not be monotone.

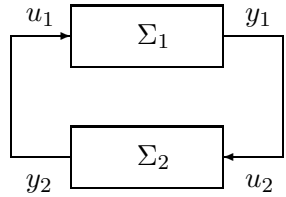

Figure 2: Systems in feedback

Based on Proposition 5.10, one of our main results will be the formulation of a small-gain theorem for the feedback interconnection of a system with monotonically increasing I/O static gain (positive path) and a system with monotonically decreasing I/O gain (negative path).

**Theorem 2** *Consider the following interconnection of two SISO dynamical systems*

$$\begin{array}{rcl} \dot{x} & = & f_x(x, w), \quad y = h_x(x) \\ \dot{z} & = & f_z(z, y), \quad w = h_z(z) \end{array} \qquad (32)$$

*with $\mathcal{U}_x = \mathcal{Y}_z$ and $\mathcal{U}_z = \mathcal{Y}_x$. Suppose that:*

1. *the first system is monotone when its input $w$ as well as output $y$ are ordered according to the standard order $\leq$ of the real axis;*

2. *the second system is monotone when its input $y$ is ordered according to $\leq$ and its output $w$ is ordered by $\geq$;*

3. *the respective static I/S characteristics $k_x(\cdot)$ and $k_z(\cdot)$ exist (thus, the static I/O characteristics $k_y(\cdot)$ and $k_w(\cdot)$ exist too and are respectively monotonically increasing and monotonically decreasing); and*

4. *every solution of the closed-loop system is bounded.*

*Then, system (32) has a globally attractive equilibrium provided that the following scalar discrete time dynamical system, evolving in $\mathcal{U}_x$:*

$$u_{k+1} = (k_w \circ k_y)(u_k) \qquad (33)$$

*has a unique globally attractive equilibrium $\bar{u}$.*

For a graphical interpretation of condition (33) see Fig. 3.

*Proof.* Equilibria of (32) are in one to one correspondence with solutions of $k_w(k_y(u)) = u$, viz. equilibria of (33). Thus, existence and uniqueness of the equilibrium follows from the GAS assumption on (33).

We need to show that such an equilibrium is globally attractive. Let $\xi \in \mathbb{R}^{n_x} \times \mathbb{R}^{n_z}$ be an arbitrary initial condition and let $y_+ := \limsup_{t \to +\infty} y(t, \xi)$ and $y_- := \liminf_{t \to +\infty} y(t, \xi)$. Then, $w_+ := \limsup_{t \to +\infty} w(t, \xi)$ and $w_- := \liminf_{t \to +\infty} w(t, \xi)$ satisfy by virtue of the second part of Proposition 5.10, applied to the $z$-subsystem:

$$k_z(y_+) \leq w_- \leq w_+ \leq k_z(y_-). \tag{34}$$

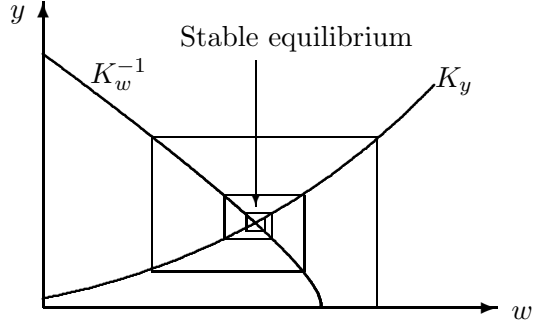

Figure 3: I/O static characteristics in the $(w, y)$ plane: negative feedback

An analogous argument, applied to the $x$-subsytem, yields: $k_x(w_-) \leq x_- \leq x_+ \leq k_x(w_+)$ and by combining this with the inequalities for $w_+$ and $w_-$ we end up with:

$$k_x(k_z(y_+)) \leq y_- \leq y_+ \leq k_x(k_z(y_-)).$$

By induction we have, after an even number $2n$ of iterations of the above argument:

$$(k_x \circ k_z)^{2n}(y_-) \leq y_- \leq y_+ \leq (k_x \circ k_z)^{2n}(y_+).$$

By letting $n \to +\infty$ and exploiting global attractivity of (33) we have $y_- = y_+$. Equation (34) yields $w_- = w_+$. Thus there exists $\bar{u}$, such that:

$$\begin{array}{rcl} \bar{u} & = & \lim_{t \to +\infty} y(t, \xi) \\ k_z(\bar{u}) & = & \lim_{t \to +\infty} w(t, \xi). \end{array} \tag{35}$$

Let $z_e$ be the (globally asymptotically stable) equilibrium (for the $z$-subsystem) corresponding to the constant input $y(t) \equiv \bar{u}$ and $x_e$ the equilibrium for the $x$-subsystem relative to the input $w(t) \equiv k_z(\bar{u})$. Clearly $\eta := [x_e, z_e]$ is the unique equilibrium of (32). The fact that $[x(t, \xi), z(t, \xi)] \to \eta$ now follows from Proposition 5.7. ∎

**Remark 6.1** We remark that traditional small-gain theorems, see for instance [16], also provide sufficient conditions for global existence and boundedness of solutions. In this respect it is of interest to notice that, for monotone systems, boundedness of trajectories follows at once provided that at least one of the interconnected systems has a uniformly bounded output map (this is always the case for instance if the state space of the corresponding system is compact). However, when both output maps are unbounded, boundedness of trajectories needs to be proved with different techniques. The following Proposition addresses this issue and provides additional conditions which together with the small-gain condition allow to conclude boundedness of trajectories. □

We say that the I/S characteristic $k_x(\cdot)$ is *unbounded* (relative to $X$) if for all $\xi \in X$ there exist $u_1, u_2 \in \mathcal{U}$ so that $k_x(u_1) \succeq \xi \succeq k_x(u_2)$.

**Lemma 6.2** Suppose that the system (1) is endowed with an unbounded I/S static characteristic and that inputs are scalar ($\mathbb{B}_\mathcal{U} = \mathbb{R}$ with the usual order). Then, for any $\xi \in X$ there exists $\bar{\xi} \in X$ so that for any input $u$:

$$\phi(t, \xi, u) \preceq \max\left\{\bar{\xi}, k_x\left(\sup_{\tau \in [0,t]} u(\tau)\right)\right\} \quad \forall\, t \geq 0. \tag{36}$$

An analogous property holds with $\preceq$ replaced by $\succeq$ and sup's replaced by inf's.

We remark that $x \preceq \max\{y_1, y_2\}$ is just a short-hand notation for $x \preceq y_1$ or $x \preceq y_2$.

*Proof.* Let $\xi \in X$ be arbitrary. As $k_x$ is unbounded there exists $\bar{u}$ such that $\xi \preceq k_x(\bar{u}) := \bar{\xi}$. Pick any input $u$ and any $t_0 \geq 0$, and let $\mu := \sup_{\tau \in [0,t_0]} u(\tau)$. There are two possibilities: $\mu \leq \bar{u}$ or $\mu \geq \bar{u}$. By monotonicity with respect to initial conditions and inputs, the first case yields:

$$\phi(t_0, \xi, u) \preceq \phi(t_0, \bar{\xi}, \bar{u}) = \bar{\xi}. \tag{37}$$

So we assume from now on that $\mu \geq \bar{u}$. We introduce the input $U$ defined as follows: $U(t) := \mu$ for all $t \leq t_0$, and $U(t) = u(t)$ for $t > t_0$. Notice that $U \succeq u$, and also that $\phi(t_0, k_x(\mu), U) = k_x(\mu)$, because the state $k_x(\mu)$ is by definition an equilibrium of $\dot{x} = f(x.\mu)$ and $U(t) \equiv \mu$ on the interval $[0, t_0]$. We conclude that

$$\phi(t_0, \xi, u) \preceq \phi(t_0, k_x(\bar{u}), U)) = k_x(\mu) \tag{38}$$

and (36) follows combining (37) and (38). The statement for $\succeq$ is proved in the same manner. ∎

**Proposition 6.3** Consider the feedback interconnection of two SISO monotone dynamical systems as in (32), and assume that the orders in both state-spaces are bounded. Assume that the systems are endowed with *unbounded* I/S static characteristics $k_x(\cdot)$ and $k_z(\cdot)$ respectively. If the small gain condition of Theorem 2 is satisfied then solutions exist for all positive times, and are bounded.

Clearly, the above result allows to apply Theorem 2 also to classes of monotone systems for which boundedness of trajectories is not apriori known.

*Proof.* We show first that solutions are upper-bounded. A symmetric argument can be used for determining a lower bound. Let $\xi, \zeta$ be arbitrary initial conditions for the $x$ and $z$ subsystems. Correspondingly solutions are maximally defined over some interval $[0, T)$. Let $t$ be arbitrary in $[0, T)$. By Lemma 6.2, equation (36) holds, for each of the systems. Moreover, composing (36) (and its counter-part for lower-bounds) with the output map yields, for suitable constants $\bar{y}, \bar{w}, \underline{y}, \underline{w}$ which only depend upon $\xi, \zeta$.

$$y(t, \xi, w) \leq \max\left\{\bar{y}, k_y\left(\max_{\tau \in [0,t]} w(\tau)\right)\right\} \tag{39}$$

$$w(t, \zeta, y) \leq \max\left\{\bar{w}, k_w\left(\min_{\tau \in [0,t]} y(\tau)\right)\right\} \tag{40}$$

$$y(t, \xi, w) \geq \min\left\{\underline{y}, k_y\left(\min_{\tau \in [0,t]} w(\tau)\right)\right\} \tag{41}$$

$$w(t, \zeta, y) \geq \min\left\{\underline{w}, k_w\left(\max_{\tau \in [0,t]} y(\tau)\right)\right\}. \tag{42}$$

Substituting equation (40) into (39) gives:

$$y(t, \xi, w) \leq \max\left\{\bar{y}, k_y(\bar{w}), k_y \circ k_w\left(\min_{\tau \in [0,t]} y(\tau)\right)\right\}. \tag{43}$$

Analogously, substitution of (41) into (43) yields (using that $k_y \circ k_w$ is a nonincreasing function):

$$y(t, \xi, w) \leq \max\left\{\bar{y}, k_y(\bar{w}), k_y \circ k_w(\underline{y}), k_y \circ k_w \circ k_y(\min_{\tau \in [0,t]} w(\tau))\right\}. \tag{44}$$

Finally equation (42) into (44) yields:

$$y(t, \xi, w) \leq \max\left\{a, \rho \circ \rho\left(\max_{\tau \in [0,t]} y(\tau)\right)\right\}, \tag{45}$$

where we are denoting $\rho := k_y \circ k_w$ and

$$a = a(\xi, \zeta) := \max\left\{\bar{y}, k_y(\bar{w}), k_y \circ k_w(\underline{y}), k_y \circ k_w \circ k_y(\underline{w})\right\}.$$

Let $y_e$ be the output value of the $x$-subsystem, corresponding to the unique equilibrium of the feedback interconnection (32).

Notice that attractivity of (33) implies attractivity of $y(t+1) = k_y \circ k_w(y(t)) := \rho(y(t))$ and a fortriori of
$$y(t+1) = \rho \circ \rho(y(t)). \tag{46}$$
We claim as follows:
$$y > y_e \quad \Rightarrow \quad \rho \circ \rho(y) < y. \tag{47}$$
By attractivity, there exists some $y_1 > y_e$ such that $\rho \circ \rho(y_1) - y_1 < 0$ (otherwise all trajectories of (46) starting from $y > y_e$ would be monotonically increasing, which is absurd). Now, assume by contradiction that there exists also some $y_2 > y_e$ such that $\rho \circ \rho(y_2) - y_2 > 0$. Then, as $\rho$ is a continuous function, there would exist an $y_0 \in (y_1, y_2)$ (or in $(y_2, y_1)$ if $y_2 < y_1$) such that $\rho \circ \rho(y_0) = y_0$. This clearly violates attractivity (at $y_e$) of (46), since $y_0$ is an equilibrium point. So (47) is proved.

Let $M := \max_{\tau \in [0,t]} y(\tau)$, so $M = y(\tau_0)$ for some $\tau_0 \in [0,t]$. Therefore (45) at $t = \tau_0$ says that $y(\tau_0) \leq \max\{a, \rho \circ \rho(y(\tau_0))\}$, and (47) applied at $y = y(\tau_0)$ gives that $y(\tau_0) \leq \max\{a, y_e\}$ (by considering separately the cases $y(\tau_0) > y_e$ and $y(\tau_0) \leq y_e$). As $y(t) \leq y(\tau_0)$, we conclude that $y(t) \leq \max\{a, y_e\}$. This shows that $y$ is upper bounded by a function which depends only on the initial states of the closed-loop system. Analogous arguments can be used in order to show that $y$ is lower bounded, and by symmetry the same applies to $w$. Thus, over the interval $[0, T)$ the $x$ and $z$ subsystems are fed by bounded inputs and by monotonicity (together with the existence of I/S static characteristics) this implies, by Proposition 5.6, that $T = +\infty$ and that trajectories are uniformly bounded. ∎

## 7 An Application

A large variety of eukaryotic cell signal transduction processes employ "Mitogen-activated protein kinase (MAPK) cascades," which play a role in some of the most fundamental processes of life (cell proliferation and growth, responses to hormones, etc). A MAPK cascade is a cascade connection of three SISO systems, each of which is (after restricting to stoichiometrically conserved subsets) either a one- or a two-dimensional system, see [12, 17]. We will show here that the two-dimensional case gives rise to monotone systems which admit static I/O characteristics. (The same holds for the much easier one-dimensional case, as follows from the results in [27].)

After nondimensionalization, the basic system to be studied is a system as in (13), where the functions $\theta_i$ are of the type
$$\theta_i(r) = \frac{a_i r}{1 + b_i r}$$
for various positive constants $i$. It follows from Proposition 3.5 that our systems (with output $y$) are monotone, and therefore *every MAPK cascade is monotone.*

We claim, further, that each such system has a static I/O characteristic. It will follow, by basic properties of cascades of stable systems, that the cascades have the same property. Thus, the complete theory developed in this paper, including small gain theorems, can be applied to MAPK cascades.

**Proposition 7.1** *For any system of the type (13), and each constant input $u$, there exists a unique globally asymptotically stable equilibrium inside $\Delta$.*

*Proof.* As the set $\Delta$ is positively invariant, the Brower Fixed-Point Theorem ensures existence of an equilibrium. In order to show that this equilibrium is unique, we consider the Jacobian of $f$:
$$\nabla f = \begin{bmatrix} -D\theta_1(x) - D\theta_2(1-x-y) & -D\theta_2(1-x-y) \\ -D\theta_3(1-x-y) & -D\theta_4(y) - D\theta_3(1-x-y) \end{bmatrix} \tag{48}$$
Notice that:
$$\operatorname{tr}(\nabla f) = -D\theta_1(x) - D\theta_2(1-x-y) - D\theta_4(y) - D\theta_3(1-x-y) < 0, \quad \forall x, y \in \Delta$$
$$\det(\nabla f) = D\theta_1(x)D\theta_4(y) + D\theta_1(x)D\theta_3(1-x-y) + D\theta_2(1-x-y)D\theta_4(y) > 0, \quad \forall x, y \in \Delta$$

The Jacobian Theorem by R. Fessler (see [10]) ensures injectivity of $f$ whenever $\nabla f$ is injective globally over $\Delta$ (as it is the case, since $\det \neq 0$). Therefore the equilibrium is unique. As $\nabla f$ is Hurwitz it is also locally asymptotically stable.

In order to conclude global asymptotic stability, we only need by Poincaré-Bendixson to rule out the possibility of limit cycles. Bendixson's Criterion shows, indeed, that no periodic orbits can exist in $\Delta$, since $\text{div}(f) < 0$. ∎

Fig. 7 shows the phase plane of the system (the diagonal line indicates the boundary of the triangular region of interest), when coefficients have been chosen so that the equations are:

$$\dot{x} = -1.0 \frac{x}{1+x} + 2 \frac{1-x-y}{3-x-y}, \quad \dot{y} = \frac{1-x-y}{2-x-y} - 2 \frac{y}{2+y}$$

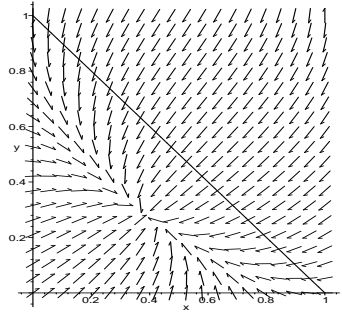

Figure 4: Direction field for example

## 8 Relations to Positivity

In this section we investigate the relationship between the notions of cooperative and positive systems. Positive linear systems (in continuous as well as discrete time) have attracted much attention in the control literature, see for instance [8, 9, 19, 21, 22, 31]. We will say that a finite dimensional linear system, possibly time-varying,

$$\dot{x} = A(t)x + B(t)u \qquad (49)$$

(where the entries of the $n \times n$ matrix $A$ and the $n \times m$ matrix $B$ are Lebesgue measurable locally essentially bounded functions of time) is *positive* if the positive orthant is forward invariant for positive input signals; in other words, for any $\xi \succeq 0$ and any $u(t) \succeq 0$ ($\succeq$ denotes here the partial orders induced by the positive orthants), and any $t_0 \in \mathbb{R}$ it holds that $\phi(t, t_0, \xi, u) \succeq 0$ for all $t \geq t_0$.

Let say that (49) is a *Metzler system* if $A(t)$ is a Metzler matrix, i.e., $A_{ij}(t) \geq 0$ for all $i \neq j$, and $B_{ij}(t) \geq 0$ for all $i, j$, for almost all $t \geq 0$. It is well known for time-invariant systems ($A$ and $B$ constant), see for instance [19], Chapter 6, or [8] for a recent reference, that a system is positive if and only if it is a Metzler system. This also holds for the general case, and we provide the proof here for completeness. For simplicity in the proof, and because we only need this case, we make a continuity assumption in one of the implications.

**Lemma 8.1** If (49) is a Metzler system then it is positive. Conversely, if (49) is positive and $A(\cdot)$ and $B(\cdot)$ are continuous, then (49) is a Metzler system.

*Proof.* Let us prove sufficiency first. Consider first any trajectory $x(\cdot)$ with $x(s) \gg 0$, any fixed $T > s$, and any input $u(\cdot)$ so that $u(t) \geq 0$ for all $t \geq s$. We need to prove that $x(T) \succeq 0$. Since $A(t)$ is essentially bounded (over any bounded time-interval) and Metzler, there is an $r > 0$ such that $rI + A(t) \geq 0$ for almost all $t \in [s, T]$, where "$\geq$" is meant elementwise. Consider $z(t) := \exp(r(t-s))x(t)$ and $v(t) := \exp(r(t-s))u(t)$, and note that $z(s) = x(s) \gg 0$ and $v(t) \geq 0$ for all $t \geq s$. We claim that $z(t) \succeq 0$ for all $t \in [s, T]$. Let $\tau > s$ be the infimum of the set of $t$'s such that $z(t) \not\succeq 0$ and assume, by contradiction, $\tau < +\infty$. By continuity of trajectories, $z(\tau) \succeq 0$. Moreover $z(\tau) = z(s) + \int_s^\tau \dot{z}(t)dt = z(s) + \int_s^\tau (rI + A(t))z(t) + B(t)v(t)\,dt \succeq z(s) \gg 0$, and therefore there exists an interval $[\tau, \tau + \varepsilon]$ such that $z(t) \succeq 0$ for all $t \in [\tau, \tau + \varepsilon]$. But this is a contradiction, unless $\tau = +\infty$ as claimed. By continuous dependence with respect to initial conditions, and closedness of the positive orthant, the result carries over to any initial condition $x(s) \succeq 0$.

For the converse implication, denote with $\Phi(t,s)$ the fundamental solution associated to $A(t)$ (viz. $\partial \Phi/\partial t = A(t)\Phi$, and $\Phi(s,s) = I$). Using $u \equiv 0$ we know that $\Phi(t,s) \geq 0$ whenever $t \geq s$ ("$\geq$" is meant here elementwise). Therefore also $[\Phi(t,s) - I]_{ij} \geq 0$ for all $i \neq j$. Since

$$A(\tau) = \left.\frac{\partial}{\partial t}\right|_{t=\tau} \Phi(t,\tau) = \lim_{t \to 0} \frac{\Phi(t,\tau) - I}{t},$$

for all $\tau$, it follows that $A(\tau)_{ij} \geq 0$ for all $i \neq j$. Consider a solution with $x(s) = 0$, $u$ constant $\geq 0$, for $t \geq s$. Since $x(t) \succeq 0$, also $(1/(t-s))x(t) \succeq 0$, and therefore, taking limits as $t \searrow 0$, $\dot{x}(s) \succeq 0$ (the derivative exists by the continuity assumption). But $\dot{x}(s) = A(s)x(s) + B(s)u$, and $x(s) = 0$, so $B(s)u \succeq 0$ for all such $u$, i.e. $B(s) \succeq 0$. ∎

Thus, by virtue of Theorem 3.3, a time-invariant linear system is cooperative if and only if it is positive. The result is a system-theoretic analog of the fact that a differentiable scalar real function is monotonically increasing if and only if its derivative is always nonnegative.

We next provide a result which says, in essence, that cooperativity is the same as "incremental positivity" (see [2], [30] for other notions of incrementality). This fact, for the special case of linear systems can be seen as a direct consequence of linearity, using the fact that the relation $\xi_1 \succeq \xi_2$ is the same as $\xi_1 - \xi_2 \succeq 0$ and similarly for inputs.

**Definition 8.2** We say that a system (1) is *incrementally positive* (or "variationally positive") if, for every solution $x(t) = \phi(t, \xi, u)$ of (1), the linearized system

$$\dot{z} = A(t)z + B(t)v \tag{50}$$

where $A(t) = \frac{\partial f}{\partial x}(x(t), u(t))$ and $B(t) = \frac{\partial f}{\partial u}(x(t), u(t))$, is a positive system. □

**Proposition 8.3** Suppose that $\mathbb{B}_{\mathcal{U}} = \mathbb{R}^m$, $\mathcal{U}$ satisfies an approximability property, and that both $\mathcal{V}$ and $\mathcal{W} = \text{int}\,\mathcal{U}$ are order-convex. Let $f(x,u)$ be continuously differentiable. Then system (1) is cooperative if and only if it is incrementally positive.

*Proof.* Under the given hypotheses, a system is cooperative iff $\frac{\partial f}{\partial x}(x,u)$ is a Metzler matrix, and every entry of $\frac{\partial f}{\partial u}(x,u)$ is nonnegative, for all $x \in X$ and all $u \in \mathcal{U}$, cf. Proposition 3.3. Therefore, by the criterion for positivity of linear time-varying systems, this implies that (50) is a positive linear time-varying system along any trajectory of (1).

Conversely, pick an arbitrary $\xi$ in $X$ and any input of the form $u(\cdot) = \bar{u} \in \mathcal{U}$. Suppose that (50) is a positive linear time-varying system along the trajectory $x(t) = x(t,\xi,u)$ (this system has continuous matrices $A$ and $B$ because $u$ is constant). Then, by the positivity criterion of linear time-varying systems, for all $t \geq 0$ we have $\frac{\partial f}{\partial x}(x(t),\bar{u})$ is Metzler and $\frac{\partial f}{\partial u}(x(t),\bar{u}) \succeq 0$. Finally, evaluating the Jacobian at $t = 0$ yields that $\frac{\partial f}{\partial x}(\xi,\bar{u})$ is Metzler and $\frac{\partial f}{\partial u}(\xi,\bar{u})$ is nonnegative. Since $\xi$ and $\bar{u}$ were arbitrary, we have the condition for cooperativity given in Proposition 3.3. ∎

**Remark 8.4** Looking at cooperativity as a notion of "incremental positivity" one can provide an alternative proof of the infinitesimal condition for cooperativity, based on the positivity of the variational equation. Indeed, assume that each system (50) is a positive linear time-varying system, along trajectories of (1). Pick arbitrary initial conditions $\xi_1 \succeq \xi_2 \in X$ and inputs $u_1 \geq u_2$. Let $\Phi(h)$ be defined as:

$$\Phi(h) := \phi(t, \xi_2 + h(\xi_1 - \xi_2), u_2 + h(u_1 - u_2)).$$

By virtue of Theorem 1 in [26] we have

$$\phi(t,\xi_1,u_1) - \phi(t,\xi_2,u_2) = \Phi(1) - \Phi(0) = \int_0^1 \Phi'(h)\,dh = \int_0^1 z_h(t,\xi_1 - \xi_2, u_1 - u_2)\,dh \tag{51}$$

where $z_h$ denotes the solution of (50) when $\frac{\partial f}{\partial u}(x,u)$ and $\frac{\partial f}{\partial u}(x,u)$ are evaluated along $\phi(t, \xi_2 + h(\xi_1 - \xi_2), u_2 + h(u_1 - u_2))$. Therefore, by positivity, exploiting (51) and monotonicity of the integral we have $\phi(t, \xi_1, u_1) - \phi(t, \xi_2, u_2) \succeq 0$, as claimed. □

We remark that monotonicity with respect to other orthants correspond to generalized positivity properties for linearizations, as should be clear by Corollary 3.4.

# Appendix: A Lemma on Invariance

We present here a characterization of invariance of relatively closed sets, under differential inclusions. The result is a simple adaptation of a well-known condition, and is expressed in terms of appropriate tangent cones.

## A.1 Differential inclusions

We let $\mathcal{V}$ be an open subset of some Euclidean space $\mathbb{R}^n$ and consider set-valued mappings $F$ defined on $\mathcal{V}$: these are mappings which assign some subset $F(x) \subseteq \mathbb{R}^n$ to each $x \in \mathcal{V}$. Associated to such mappings $F$ are *differential inclusions*

$$\dot{x} \in F(x) \qquad (52)$$

and one says that a function $x : [0, T] \to \mathcal{V}$ is a *solution* of (52) if $x$ is an absolutely continuous function with the property that $\dot{x}(t) \in F(x(t))$ for almost all $t \in [0, T]$. A set-valued mapping $F$ is *compact-valued* if $F(\xi)$ is a compact set, for each $\xi \in \mathcal{V}$.

A set-valued mapping $F$ is said to be *locally Lipschitz* if the following property holds: for each compact subset $C \subseteq \mathcal{V}$ there is some constant $k$ such that

$$F(\xi) \subseteq F(\zeta) + k\, |\xi - \zeta|\, B \qquad \forall \xi, \zeta \in C$$

where $B$ denotes the unit ball in $\mathbb{R}^n$. In other words, for any two elements $\xi, \zeta$ of $C$, and for every $f \in F(\xi)$, there must exist some $f' \in F(\zeta)$ such that $|f - f'| \le k\, |\xi - \zeta|$. (We use $|x|$ to denote Euclidean norm in $\mathbb{R}^n$.) Note that when $F(x) = \{f(x)\}$ is single-valued, this is the usual definition of a locally Lipschitz function. More generally, suppose that $f(x, u)$ is locally Lipschitz in $x \in \mathcal{V}$, locally uniformly on $u$, and pick any compact subset $D$ of the input set $\mathcal{U}$; then $F_D(x) = \{f(x, u), u \in D\}$ is locally Lipschitz. Indeed, given any compact $C \subseteq \mathcal{V}$, we may pick a constant $k$ such that $|f(\xi, u) - f(\zeta, u)| \le k\, |\xi - \zeta|$ for all $\xi, \zeta \in C$ and all $u \in D$ (local Lipschitz condition on $f$), and thus, given $f = f(\xi, u)$, we verify the definition by picking $f' := f(\zeta, u)$ with the same $u$. Also, $F_D$ is compact-valued, since each set $f(\xi, D)$ is compact, as $f$ is continuous.

We say that the set-valued mapping $F$ defined on $\mathcal{V}$ is *locally bounded* if for each compact subset $C \subseteq \mathcal{V}$ there is some constant $k$ such that

$$F(\xi) \subseteq k\, B \qquad \forall \xi \in C\,.$$

When $F$ has the form $F_D$ as above, it is locally bounded, since $F_D(\xi) \subseteq f(C \times D)$, and, $f$ being continuous, the latter set is compact.

## A.2 Invariance

Let $\mathcal{S}$ be a (nonempty) closed subset relative to $\mathcal{V}$, that is, $\mathcal{S} = S \bigcap \mathcal{V}$ for some closed subset $S$ of $\mathbb{R}^n$. We wish to characterize the property that solutions which start in the set $\mathcal{S}$ must remain there:

**Definition A.5** The subset $\mathcal{S}$ is *strongly invariant* under the differential inclusion (52) if the following property holds: for every solution $x : [0, T] \to \mathcal{V}$ which has the property that $x(0) \in \mathcal{S}$, it must be the case that $x(t) \in \mathcal{S}$ for all $t \in [0, T]$. □

For any $\xi \in \mathcal{S}$, the tangent cone to $\mathcal{S}$ at $\xi$ was introduced in Definition 3.2. Note that a vector $v$ belongs to $\mathcal{T}_\xi \mathcal{S}$ if and only if there is a sequence of elements $v_i \in \mathcal{V}$, $v_i \to v$ and a sequence $t_i \searrow 0$ such that $x + t_i v_i \in \mathcal{S}$ for all $i$. This is a standard definition, see e.g. [5], Section 2.7, and the object so defined is usually called the *Bouligand* or *contingent* tangent cone (to differentiate it from other types of tangent cones). The set $\mathcal{T}_\xi\mathcal{S}$ is indeed a cone ($v \in \mathcal{T}_\xi\mathcal{S}$ and $\lambda \in \mathbb{R}_{\geq 0}$ imply $\lambda v \in \mathcal{T}_\xi\mathcal{S}$: just use $t_i/\lambda$ if $\lambda > 0$, and $\xi_i \equiv \xi$ if $\lambda = 0$) and is always closed (if $v_j \to v$ and $(1/t_{ij})\left(\xi_{ij} - \xi\right) \to v_j$ with $\xi_{ij} \underset{\mathcal{S}}{\to} \xi$ and $t_{ij} \searrow 0$ as $i \to \infty$, then one may pick a sequence $\{i(j)\}$ such that $\left|(1/t_{i(j)j})\left(\xi_{i(j)j} - \xi\right) - v_j\right| < 1/j$ and $t_{i(j)j} \searrow 0$ as $j \to \infty$; this sequence satisfies that $(1/t_{i(j)j})\left(\xi_{i(j)j} - \xi\right) \to v$). Further, $\mathcal{T}_\xi\mathcal{S} = \mathbb{R}^n$ when $x$ is in the interior of $\mathcal{S}$ relative to $\mathcal{V}$ (so only boundary points are of interest).

**Theorem 3** *Suppose that $F$ is a locally Lipschitz, compact-valued, and locally bounded set-valued mapping on the open subset $\mathcal{V} \subseteq \mathbb{R}^n$, and $\mathcal{S}$ is a closed subset of $\mathcal{V}$. Then, the following two properties are equivalent:*

1. *$\mathcal{S}$ is strongly invariant under $F$.*
2. *$F(\xi) \subseteq \mathcal{T}_\xi\mathcal{S}$ for every $\xi \in \mathcal{S}$.*

Just for purposes of the proof, let us say that a set-valued mapping $F$ is "nice" if $F$ is defined on all of $\mathbb{R}^n$ and it satisfies the following properties: $F$ is locally Lipschitz, compact-valued, convex-valued, and globally bounded ($F(\xi) \subseteq kB$ for all $\xi \in \mathbb{R}^n$, for some $k$). Theorem 4.3.8 in [5] establishes that Properties 1 and 2 in the statement of Theorem 3 are equivalent, and are also equivalent to:

$$F(\xi) \subseteq \operatorname{co}\mathcal{T}_\xi\mathcal{S} \quad \text{for every } \xi \in \mathcal{S} \tag{53}$$

("co" indicates closed convex hull) provided that $\mathcal{S}$ is a closed subset of $\mathbb{R}^n$ and $F$ is nice (a weaker linear growth condition can be replaced for global boundedness, c.f. the "standing hypotheses" in Section 4.1.2 of [5]). We will reduce to this case using the following observation.

**Lemma A.6** *Suppose that $F$ is a locally Lipschitz, compact-valued, and locally bounded set-valued mapping on the open subset $\mathcal{V} \subseteq \mathbb{R}^n$, and $\mathcal{S}$ is a closed subset of $\mathcal{V}$. Let $M$ be any given compact subset of $\mathcal{V}$. Then, there exist a nice set-valued $\widehat{F}$ and a closed subset $\mathcal{S}'$ of $\mathbb{R}^n$ such that the following properties hold:*

$$F(\xi) \subseteq \widehat{F}(\xi) \qquad \forall \xi \in M \tag{54}$$

$$M \bigcap \mathcal{S} \subseteq \mathcal{S}' \subseteq \mathcal{S} \tag{55}$$

$$\forall \xi \in \mathcal{S}', \text{ either } \widehat{F}(\xi) = \{0\} \text{ or } \mathcal{T}_\xi\mathcal{S} = \mathcal{T}_\xi\mathcal{S}' \tag{56}$$

$$\mathcal{T}_\xi\mathcal{S} = \mathcal{T}_\xi\mathcal{S}' \qquad \forall \xi \in M. \tag{57}$$

*and*

$$\mathcal{S} \text{ strongly invariant under } F \quad \Rightarrow \quad \mathcal{S}' \text{ strongly invariant under } \widehat{F} \tag{58}$$

*Proof.* Consider the convexification $\widetilde{F}$ of $F$; this is the set-valued function on $\mathcal{V}$ which is obtained by taking the convex hull of the sets $F(\xi)$, i.e. $\widetilde{F}(\xi) := \operatorname{co} F(\xi)$ for each $\xi \in \mathcal{V}$. It is an easy exercise to verify that if $F$ is compact-valued, locally Lipschitz, and locally bounded, then $\widetilde{F}$ also has these properties.

Clearly, if $\mathcal{S}$ is strongly invariant under $\widetilde{F}$ then it is also strongly invariant under $F$, because every solution of $\dot{x} \in F(x)$ must also be a solution of $\dot{x} \in \widetilde{F}(x)$. Conversely, suppose that $\mathcal{S}$ is strongly invariant under $F$, and consider any solution $x : [0, T] \to \mathcal{V}$ of $\dot{x} \in \widetilde{F}(x)$ which has the property that $x(0) \in \mathcal{S}$. The Filippov-Ważewski Relaxation Theorem (see e.g. [4, 11]) provides a sequence of solutions $x_k$, $k = 1, 2, \ldots$, of $\dot{x} \in F(x)$ on the interval $[0, T]$, with the property that $x_k(t) \to x(t)$ uniformly on $t \in [0, T]$ and also $x_k(0) = x(0) \in \mathcal{S}$ for all $k$. Since $\mathcal{S}$ is strongly invariant under $F$, it follows that $x_k(t) \in \mathcal{S}$ for all $k$ and $t \in [0, T]$, and taking the limit as $k \to \infty$ this implies that also $x(t) \in \mathcal{S}$ for all $t$. In summary, invariance under $F$ or $\widetilde{F}$ are equivalent, for closed sets.

Let $N$ be a compact subset of $\mathcal{V}$ which contains $M$ in its interior $\operatorname{int} N$ and pick any smooth function $\varphi : \mathbb{R}^n \to \mathbb{R}_{\geq 0}$ with support equal to $N$ (that is, $\varphi(\xi) \equiv 0$ if $x \notin \operatorname{int} N$ and $\varphi(\xi) > 0$ on $\operatorname{int} N$) and such that $\varphi(\xi) \equiv 1$ on the set $M$. Now consider the new differential inclusion defined on all of $\mathbb{R}^n$ given by $\widehat{F}(\xi) := \varphi(\xi)\widetilde{F}(\xi)$ if $\xi \in N$ and equal to $\{0\}$ outside $N$. Since $\widetilde{F}$ is locally Lipschitz and locally bounded, it follows by a standard argument that $\widehat{F}$ has these same properties. Moreover, $\widehat{F}$ is globally bounded and it is also convex-valued and compact-valued (see e.g. [18]). Thus $\widehat{F}$ is nice, as required.

Note that Property (54) holds, because $F(\xi) \subseteq \widetilde{F}(\xi)$ and $\varphi \equiv 1$ on $M$.

Let $\mathcal{S}' := \mathcal{S} \bigcap N$ (cf. Figure 5); this is a closed subset of $\mathbb{R}^n$ because the compact set $N$ has a strictly positive distance to the complement of $\mathcal{V}$. Property (55) holds as well, because $M \subseteq N$.

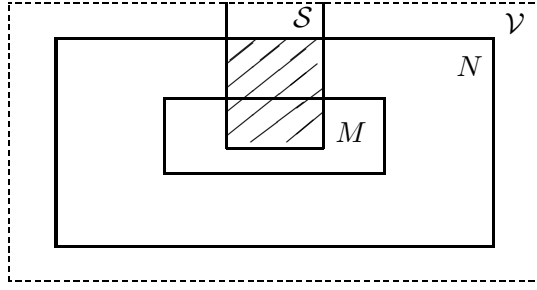

Figure 5: Shaded area is set $\mathcal{S}'$

Now pick any $\xi \in \mathcal{S}'$. There are two cases to consider: $\xi$ is in the boundary of $N$ or in the interior of $N$. If $\xi \in \partial N$, then $\widehat{F}(\xi) = \{0\}$ because $\varphi(\xi) = 0$. If instead $\xi$ belongs to the interior of $N$, there is some open subset $V \subseteq N$ such that $\xi \in V$. Therefore any sequence $\xi_i \to \xi$ with all $\xi_i \in \mathcal{S}$ has, without loss of generality, $\xi_i \in V \bigcap \mathcal{S} \subseteq N \bigcap \mathcal{S} = \mathcal{S}'$, so also $\xi_i \to \xi$ in $\mathcal{S}'$; this proves that $\mathcal{T}_\xi \mathcal{S} \subseteq \mathcal{T}_\xi \mathcal{S}'$, and the reverse inclusion is true because $\mathcal{S}' \subseteq \mathcal{S}$. Hence Property (56) has been established.

Regarding Property (57), this follows from the discussion in the previous paragraph, since $M$ is included in the interior of $N$.

In order to prove the last property in the theorem, we start by remarking that if $x : [0, T] \to \mathbb{R}^n$ is a solution of $\dot{x} \in \widehat{F}(x)$ with the property that $x(t)$ belongs to the interior of $N$ for all $t$ (equivalently, $\varphi(x(t)) \neq 0$ for all $t$), then there is a reparametrization of time such that $x$ is a solution of $\dot{x} \in \widetilde{F}(x)$. In precise terms: there is an interval $[0, R]$, an absolutely continuous function $\alpha : [0, \infty) \to [0, \infty)$ such that $\alpha(0) = 0$ and $\alpha(R) = T$, and a solution $z : [0, R] \to \mathbb{R}^n$ of $\dot{z} \in \widetilde{F}(z)$ such that $z(r) = x(\alpha(r))$ for all $r \in [0, R]$. To see this, it is enough (chain rule, remembering that $\widehat{F}(\xi) = \varphi(\xi)\widetilde{F}(x)$) for $\alpha$ to solve the initial value problem $d\alpha/dr = \beta(\alpha(r))$, $\alpha(0) = 0$, where $\beta(t) = 1/\varphi(x(t))$ for $t \leq T$ and $\beta(t) \equiv \beta(T)$ for $t > T$. The function $\varphi(x(t))$ is absolutely continuous, and is bounded away from zero for all $t \leq T$ (because the solution $x$ lies in a compact subset of the interior of the support of $\varphi$), so $\beta$ is locally Lipschitz and a (unique) solution exists. Since $\beta$ is globally bounded, the solution has no finite escape times. In addition, since the vector field is everywhere positive, $\alpha(r) \to \infty$ as $s \to \infty$, so there is some $R$ such that $\alpha(R) = T$.

Now suppose that $\mathcal{S}$ is invariant under $F$. As remarked earlier, this implies that $\mathcal{S}$ is also invariant under its convexification $\widetilde{F}$.

Suppose that $x : [0, T] \to \mathbb{R}^n$ is a solution of $\dot{x} = \widehat{F}(x)$ such that $x(0) \in \mathcal{S}'$ and $x(t)$ is in the interior of $N$ for all $t$. We find a solution $z$ of $\dot{z} \in \widetilde{F}(z)$ such that $z(r) = x(\alpha(r))$ for all $r \in [0, R]$ and $z(0) = x(0) \in \mathcal{S}' \subseteq \mathcal{S}$ as earlier. Invariance of $\mathcal{S}$ under $\widetilde{F}$ gives that $z(r)$, and hence $x(t)$, remains in $\mathcal{S}$. Since $\mathcal{S}' = \mathcal{S} \bigcap N$, we conclude that $x(t) \in \mathcal{S}'$ for all $t \in [0, T]$.

Next, we use some ideas from the proof of Theorem 4.3.8 in [5]. Pick any $\xi_0 \in \mathcal{S}'$, and any $v \in \widehat{F}(\xi_0)$. Define the mapping $f : \mathbb{R}^n \to \mathbb{R}^n$ by the following rule: for each $\xi \in \mathbb{R}^n$, $f(\xi)$ is the unique closest point to $v$ in $\widehat{F}(\xi)$. As in the above citation, this map is continuous. We claim that, for each $\xi \in \mathcal{S}'$ there is some $\delta > 0$ and a solution of $\dot{x} = f(x)$ such that $x(0) = \xi$ and $x(t) \in \mathcal{S}'$ for all $t \in [0, \delta]$. (Note that, in particular, this $x$ solves $\dot{x} \in \widehat{F}(x)$.) If $\xi$ is on the boundary of $N$, then $\widehat{F}(\xi) = \{0\}$ implies that $f(\xi) = 0$, and hence $x(t) \equiv \xi$ is such a solution. If instead $\xi$ belongs to the interior of $N$ then the previous remarks shows that $x(t) \in \mathcal{S}'$ for all $t \in [0, \delta]$, where we pick a smaller $\delta$ if needed in order to insure that $x(t)$ remains in the interior of $N$. We conclude from the claim that the closed set $\mathcal{S}'$ is locally-in-time invariant with respect to the differential inclusion $\{f(x)\}$, which satisfies the "standing hypotheses" in Chapter 4 of [5]. This inclusion is hence also "weakly invariant" as follows from Exercise 4.2.1 in that textbook. This in turn implies, by Theorem 4.2.10 there, that $\langle f(\xi), \zeta \rangle \leq 0$ for all $\xi \in \mathcal{S}'$ and all $\zeta$ in the proximal normal set $N_\xi \mathcal{S}'$ defined in that reference (we are using a different notation). Applied in particular at the point $\xi_0$ (so that $f(\xi_0) = v$), we conclude that $\langle v, \zeta \rangle \leq 0$ for all $\zeta \in N_{\xi_0} \mathcal{S}'$. Since $v$ was an arbitrary element of $\widehat{F}(\xi_0)$, it follows that the upper Hamiltonian condition in part (d) of Theorem 4.3.8 in [5] holds for the map $\widehat{F}$ at the point $\xi_0$. Since $\xi_0$ was itself an arbitrary point in $\mathcal{S}'$, the condition holds on all of $\mathcal{S}'$. Therefore $\mathcal{S}'$ is invariant for $\widehat{F}$, as claimed. ∎

## A.3  Proof of Theorem 3

**Proof that 2⇒1**

Suppose that $F(\xi) \subseteq \mathcal{T}_\xi \mathcal{S}$ for every $\xi \in \mathcal{S}$, and pick any solution $x : [0, T] \to \mathcal{V}$ of $\dot{x} \in F(x)$ with $x(0) \in \mathcal{S}$.

Since $x(\cdot)$ is continuous, there is some compact subset $M \subseteq \mathcal{V}$ such that $x(t) \in M$ for all $t \in [0, T]$. We apply Lemma A.6 to obtain $\widehat{F}$ and $\mathcal{S}'$. By Property (54), it holds that $x$ is also a solution of $\dot{x} \in \widehat{F}(x)$, and Property (55) gives that $x(0)$ belongs to the subset $\mathcal{S}'$.

Taking convex hulls, $\widetilde{F}(x) \subseteq \operatorname{co} \mathcal{T}_\xi \mathcal{S}$ for every $x \in \mathcal{S}$. Since $\widehat{F}$ is a scalar multiple of $\widetilde{F}$, and $\operatorname{co} \mathcal{T}_\xi \mathcal{S}$ is a cone (because $\mathcal{T}_\xi \mathcal{S}$ is a cone), it follows that $\widehat{F}(\xi) \subseteq \operatorname{co} \mathcal{T}_\xi \mathcal{S}$ for every $\xi \in \mathcal{S}$, and so also for $\xi \in \mathcal{S}'$. By Property (56),

$$\widehat{F}(\xi) \subseteq \operatorname{co} \mathcal{T}_\xi \mathcal{S}' \qquad \forall \xi \in \mathcal{S}', \tag{59}$$

since either $\widehat{F}(\xi) = 0$ or $\mathcal{T}_\xi \mathcal{S}' = \mathcal{T}_\xi \mathcal{S}$ (and hence their convex hulls coincide).

In summary, Property (53) is valid for $\widehat{F}$ in place of $F$ and $\mathcal{S}'$ in place of $\mathcal{S}$, and $\widehat{F}$ is nice. Thus we may apply Theorem 4.3.8 in [5] to conclude that $\mathcal{S}'$ is strongly invariant under $\widehat{F}$. Since $x(0) \in \mathcal{S}'$, it follows that $x(t) \in \mathcal{S}'$ for all $t \in [0, T]$, and therefore also $x(t) \in \mathcal{S}$ for all $t \in [0, T]$, as wanted.

**Proof that 1⇒2**

Suppose that $\mathcal{S}$ is strongly invariant under $F$, and pick any $\xi_0 \in \mathcal{S}$. We apply Lemma A.6, with $M = \{\xi_0\}$, to obtain $\widehat{F}$ and $\mathcal{S}'$. Note that $M \cap \mathcal{S} = \{\xi_0\}$, so $\xi_0 \in \mathcal{S}'$. Then Property (58) gives that $\mathcal{S}'$ is strongly invariant under $\widehat{F}$. Since $\mathcal{S}'$ is closed and $\widehat{F}$ is nice, Theorem 4.3.8 in [5] gives that $\widehat{F}(\xi) \subseteq \mathcal{T}_\xi \mathcal{S}'$ for all $\xi \in \mathcal{S}'$, and in particular for $\xi = \xi_0$. By Property (56), either $\widehat{F}(\xi_0) = \{0\}$ or $\mathcal{T}_\xi \mathcal{S}' = \mathcal{T}_\xi \mathcal{S}$, so we have that $\widehat{F}(\xi) \subseteq \mathcal{T}_\xi \mathcal{S}$ for $\xi = \xi_0$. Moreover, Property (54) gives that $F(\xi) \subseteq \widehat{F}(\xi)$ for $\xi = \xi_0$. Since $\xi_0$ was an arbitrary element of $\mathcal{S}$, the proof is complete. ∎